\documentclass{amsart}

\usepackage{amsmath,amscd,amssymb,amsthm}
\usepackage{enumerate,mathrsfs}
\usepackage{geometry}
\usepackage{hyperref}
\usepackage{url}
\usepackage{ifthen}
\usepackage[all]{xy} \xyoption{2cell} \UseAllTwocells
\usepackage[mathscr]{euscript}
\usepackage[T1]{fontenc}
\usepackage[utf8]{inputenc}   
\usepackage{refcount}
\usepackage{zref-user}
\usepackage{zref-abspage}

\numberwithin{equation}{section}

\theoremstyle{plain}
\newtheorem{thm_}[equation]{Theorem}
\newtheorem{lemma_}[equation]{Lemma}
\newtheorem{prop_}[equation]{Proposition}
\newtheorem{cor_}[equation]{Corollary}
\newtheorem{eg_}[equation]{Example}

\theoremstyle{definition}
\newtheorem{thmu_}[equation]{Theorem}
\newtheorem*{thmus_}{Theorem}
\newtheorem{propu_}[equation]{Proposition}
\newtheorem*{propus_}{Proposition}
\newtheorem{coru_}[equation]{Corollary}
\newtheorem*{corus_}{Corollary}
\newtheorem{lemu_}[equation]{Lemma}
\newtheorem*{lemus_}{Lemma}
\newtheorem{egu_}[equation]{Example}
\newtheorem*{egus_}{Example}
\newtheorem{def_}[equation]{Definition}
\newtheorem*{defs_}{Definition}
\newtheorem{rk_}[equation]{Remark}
\newtheorem*{rks_}{Remark}
\newtheorem{ex_}[equation]{Remark}
\newtheorem*{exs_}{Remark}
\newtheorem{constr_}[equation]{Construction}
\newtheorem*{constrs_}{Construction}
\newtheorem{nota_}[equation]{Notation}
\newtheorem*{notas_}{Notation}

\newcommand{\thm}[1]{\begin{thm_}#1\end{thm_}}
\newcommand{\thmu}[1]{\begin{thmu_}#1\end{thmu_}}

\newcommand{\lemm}[1]{\begin{lemma_}#1\end{lemma_}}
\newcommand{\lemu}[1]{\begin{lemu_}#1\end{lemu_}}

\newcommand{\prop}[1]{\begin{prop_}#1\end{prop_}}

\newcommand{\defi}[1]{\begin{def_}#1\end{def_}}

\newcommand{\rk}[1]{\begin{rk_}#1\end{rk_}}

\newcommand{\cor}[1]{\begin{cor_}#1\end{cor_}}

\newcommand{\constr}[1]{\begin{constr_}#1\end{constr_}}

\newcommand{\nota}[1]{\begin{nota_}#1\end{nota_}}

\newcommand{\pf}[1]{\begin{proof}#1\end{proof}}

\DeclareMathOperator{\PGL}{PGL}
\DeclareMathOperator{\SL}{SL}

\DeclareMathOperator{\Hom}{Hom}

\DeclareMathOperator{\Spec}{Spec}

\DeclareMathOperator{\im}{im}
\DeclareMathOperator{\Ob}{Ob}

\DeclareMathOperator{\Br}{Br}

\DeclareMathOperator{\res}{res}

\DeclareMathOperator{\ch}{char}

\DeclareMathOperator{\id}{id}

\DeclareMathOperator{\Fun}{Fun}

\DeclareMathOperator{\N}{N}
\DeclareMathOperator{\Map}{Map}
\newcommand{\Set}{{\sS\tu{et}}}
\newcommand{\Cat}{{\sC\tu{at}}}
\newcommand{\sPr}{{\sP\tu{r}}}
\newcommand{\PrL}{{\sPr^\tu{L}}}
\newcommand{\PrR}{{\sPr^\tu{R}}}
\newcommand{\Sch}{{\sS\tu{ch}}}
\newcommand{\Esp}{{\sE\tu{sp}}}
\newcommand{\Chp}{{\sC\tu{hp}}}
\newcommand{\bDel}{{\mathbf{\Delta}}}
\newcommand{\tril}{\triangleleft}

\newcommand{\RAd}{{\tu{RAd}}}

\newcommand{\ZZ}{\mathbb Z}
\newcommand{\bfA}{\mathbf A}%
\newcommand{\bfG}{\mathbf G}%

\newcommand{\sA}{\mathscr A}
\newcommand{\sB}{\mathscr B}
\newcommand{\sC}{\mathscr C}
\newcommand{\sD}{\mathscr D}
\newcommand{\sE}{\mathscr E}
\newcommand{\sF}{\mathscr F}

\newcommand{\sH}{\mathscr H}

\newcommand{\sK}{\mathscr K}

\newcommand{\sP}{\mathscr P}
\newcommand{\sR}{\mathscr R}
\newcommand{\sS}{\mathscr S}
\newcommand{\sT}{\mathscr T}

\newcommand{\sX}{\mathscr X}
\newcommand{\sY}{\mathscr Y}

\newcommand{\cX}{\mathcal X}

\newcommand{\s}{\sigma}

\newcommand{\La}{\Lambda}
\newcommand{\lam}{\lambda}%
\newcommand{\ep}{\epsilon}%
\newcommand{\tm}{\times}%

\newcommand{\ol}{\overline}

\newcommand{\wt}{\widetilde}
\newcommand{\bu}{\bullet}

\newcommand{\ra}{\rightarrow}

\newcommand{\xra}{\xrightarrow}

\newcommand{\hra}{\hookrightarrow}

\newcommand{\mpt}{\mapsto}

\newcommand{\is}[2]{\xymatrix@-4mm{#1 \ar[r]^-{\sim} & #2 }}
\newcommand{\mis}[2]{\xymatrix@-2mm{#1 \ar[r]^-{\sim} & #2 }}

\newcommand{\dra}[4]{\xymatrix@-4mm{#1 \ar@<.5ex>[r]^-{#3} \ar@<-.5ex>[r]_-{#4}& #2 }}  
\newcommand{\era}[5]{\xymatrix@-4mm{#1 \ar[r] &#2 \ar@<.5ex>[r]^-{#4} \ar@<-.5ex>[r]_-{#5}& #3 }}  
\newcommand{\ds}{{/\mkern-3mu/}} 

\newcommand{\tc}[1]{{\textcircled{#1}}}

\newcommand{\tu}[1]{\text{\upshape #1}}

\newcommand{\op}{{op}}

\newcommand{\et}{\tu{\'et}}
\newcommand{\liset}{\tu{lis-\'et}}
\newcommand{\cart}{\tu{cart}}

\newcommand{\fppf}{\tu{fppf}}

\newcommand{\tor}{\tu{tor}}

\DeclareFontFamily{U}{wncy}{}
\DeclareFontShape{U}{wncy}{m}{n}{%
   <5>wncyr5%
   <6>wncyr6%
   <7>wncyr7%
   <8>wncyr8%
   <9>wncyr9%
   <10>wncyr10%
   <11>wncyr10%
   <12>wncyr6%
   <14>wncyr7%
   <17>wncyr8%
   <20>wncyr10%
   <25>wncyr10}{}
\DeclareMathAlphabet{\cyrille}{U}{wncy}{m}{n}
\def\Sha{\cyrille X}

\newcommand{\eq}[1]{\begin{equation}#1\end{equation}}
\newcommand{\eqn}[1]{\begin{equation*}#1\end{equation*}}
\newcommand{\ga}[1]{\begin{gather}#1\end{gather}}
\newcommand{\gan}[1]{\begin{gather*}#1\end{gather*}}

\newcommand{\enmt}[1]{\begin{enumerate}#1\end{enumerate}}

\newcommand{\aci}[1]{\ar@{^(->}[#1]|-{/}}
\newcommand{\coaci}[1]{\ar@{_(->}[#1]|-{/}}
\newcommand{\aoi}[1]{\ar@{^(->}[#1]|-{\circ}}
\newcommand{\coaoi}[1]{\ar@{_(->}[#1]|-{\circ}}
\makeatletter
\def\citet@url@sp{https://stacks.math.columbia.edu/}
\def\citet@bib@sp{stacks-project}
\def\citet@url@kd{https://kerodon.net/}
\def\citet@bib@kd{kerodon}
\newcommand{\citet@tag}[2]{\href{#2tag/#1}{#1}}
\newcommand{\citet@taglist}[2]{%
 \def\@citet@e{}%
 \def\@citet@tag@n{0}
 \@for\@citet@tag:=#1\do{%
  \edef\@citet@tag@n{\the\numexpr\@citet@tag@n + 1}%
 }%
 \def\@citet@tags{%
  \def\@citet@tag@i{0}%
  \@for\@citet@tag:=#1\do{%
   \edef\@citet@tag@i{\the\numexpr\@citet@tag@i + 1}%
   \ifthenelse{\@citet@tag@i > 1}{
    \ifthenelse{\@citet@tag@i = \@citet@tag@n}{
     \citet@seplast%
    }{%
     \citet@sep%
    }%
   }{}%
   \citet@entry{\@citet@tag}{#2}
  }%
 }%
 \ifthenelse{\@citet@tag@n > 1}{
  \def\@citet@Tag{Tags}%
 }{%
  \def\@citet@Tag{Tag}%
 }%
 \@citet@Tag~\@citet@tags
}
\newcommand{\citet@sep}{, }
\newcommand{\citet@seplast}{ and }
\newcommand{\citet@entry}[2]{\citet@tag{#1}{#2}}
\let\@old@cite\cite
\renewcommand{\cite}[2][]{%
 \def\@citet@detail{\citet@taglist{#1}{\@citet@url}}%
 \ifthenelse{\equal{#2}{sp}}{%
  \def\@citet@url{\citet@url@sp}%
  \def\@citet@bib{\citet@bib@sp}%
 }{\ifthenelse{\equal{#2}{kd}}{%
  \def\@citet@url{\citet@url@kd}%
  \def\@citet@bib{\citet@bib@kd}%
 }{
  \def\@citet@detail{#1}%
  \def\@citet@bib{#2}%
 }}%
 \ifthenelse{\equal{#1}{}}{%
  \@old@cite{\@citet@bib}%
 }{%
  \@old@cite[\@citet@detail]{\@citet@bib}%
 }%
}
\makeatother

\newcommand{\etale}{{\'etale}}

\newcommand{\Cech}{{\v{C}ech}}

\newcommand{\Grot}{{Grothendieck}}

\newcommand{\Poin}{{Poincar\'e}}

\newcommand{\Kunn}{{K\"unneth}}

\newcommand{\BM}{{Brauer--Manin}}

\newcommand{\TS}{{Tate--Shafarevich}}

\newcommand{\Cart}{{Cartesian}}

\newcommand{\adele}{{ad\`ele}}
\newcommand{\adelic}{{ad\`elic}}

\newcommand{\obs}{\tu{obs}}

\newcommand{\desc}{\tu{desc}}
\newcommand{\conn}{\tu{conn}}
\newcommand{\fin}{\tu{fin}}
\newcommand{\fdesc}{{\fin, \desc}}
\newcommand{\ddesc}{{\desc, \desc}}

\newcommand{\etBr}{{\et, {\Br}}}
\newcommand{\sdesc}{{2\tu{-}\desc}}

\newcommand{\hdesc}[2]{{#1\tu{-}\desc{\ifthenelse{\equal{#2}{}}{}{_#2}}}}

\newcommand{\XA}{X(\bfA_k)}
\newcommand{\YA}{Y(\bfA_k)}
\newcommand{\Xk}{X(k)}

\newcounter{arrowcounter}

\makeatletter
\zref@newprop{arrownum}{\thearrowcounter}
\zref@addprops{main}{arrownum}

\newcommand{\arrowlabel}[1]{%
    \refstepcounter{arrowcounter}%
      \zref@labelbyprops{#1}{arrownum}%
    \hypertarget{#1}{\tc{\thearrowcounter}}%
          }

\newcommand{\arrowref}[1]{%
            \hyperlink{#1}{\tc{\zref@extract{#1}{arrownum}}}%
              }
\makeatother

\DeclareMathOperator{\Sh}{Sh}

\newcommand{\PTop}{{\sP\sT\tu{opos}}}
\newcommand{\Rind}{{\sR\tu{ind}}}
\newcommand{\Ring}{{\sR\tu{ing}}}

\newcommand{\sr}{{\tu{sm},\tu{rep}}}
\newcommand{\ssr}{{\tu{sm},\tu{surj},\tu{rep}}}

\newcommand{\rep}{\tu{rep}}
\newcommand{\sFo}{{\sF^\obs}}
\newcommand{\sFoB}{{\sFo, \Box}}

\newcommand{\sFm}{{\sF^m}}

\newcommand{\olobs}{{\ol\obs}}
\newcommand{\sFolo}{{\sF^\olobs}}
\newcommand{\obsol}{{\obs,\olobs}}
\newcommand{\obsm}{{\obs,m}}
\newcommand{\tobs}{{\wt{\obs}}}
\newcommand{\tBr}{{\wt{\Br}}}
\newcommand{\tdesc}{{\wt{\desc}}}
\newcommand{\tconn}{{\wt{\conn}}}
\newcommand{\tsdesc}{{\wt{\sdesc}}}
\newcommand{\tho}{{\wt{h}}}
\newcommand{\thoZ}{{\wt{h\ZZ}}}
\newcommand{\tfdesc}{{\wt{\fdesc}}}
\newcommand{\tddesc}{{\wt{\ddesc}}}
\newcommand{\tetBr}{{\wt{\etBr}}}

\begin{document}
\title[Cohomological descent for obstructions]
 {Cohomological descent for obstructions to local-global principle}
\author[C. Lv]{Chang Lv}
\address{State Key Laboratory of Cyberspace Security Defense\\
Institute of Information Engineering\\
Chinese Academy of Sciences\\
Beijing 100093, P.R. China}
\email{lvchang@amss.ac.cn}

\date{\today}
\thanks{The author is partially supported by
 National Natural Science Foundation of China NSFC Grant No. 11701552.}
\begin{abstract}
We develop a formalism of
 cohomological descent encoding {\adelic} points and
 obstructions to local-global principle on algebraic stacks.
As an application, by constructing new obstructions
 using the formalism,
  we obtain some comparison results of obstructions
  on some classes of algebraic stacks.
\end{abstract}
\maketitle

\setcounter{tocdepth}{1}
\tableofcontents

\section{Introduction} \label{intro}
Let $k$ be a number field and $X$ a $k$-variety.
Then one considers local-global principle for rational points of $X$.
The most developed method is to consider obstructions to it.
Namely, people construct various ``nice" subsets $\XA^{\obs}$ of
 {\adelic} points $\XA$, such that
 $\Xk\subseteq \XA^{\obs}\subseteq \XA$ and hope that
 they ``accounts" the local-global principle.

Among others, we are interested in relations between obstructions.
One of the results is about the well-known relations
\eqn{
\XA^\desc\subseteq \XA^{\PGL} = \XA^{\Br}
}
 for regular, quasi-projective $k$-variety $X$
 (see, e.g., \cite[Prop. 8.5.3]{poonen17rational}).
Harari \cite{harari02groupes} showed that
 for  smooth geometrically integral variety,
\eqn{
\XA^{\Br} = \XA^{\conn} = \XA^{\sdesc}.
}
By works of Stoll, Skorobogatov, Demarche, Poonen, Xu and Cao,
 it is known that for $X$ being a
 smooth, quasi-projective, geometrically integral
 $k$-variety,  $\XA^\desc=\XA^\etBr$ \cite{stoll07finite,
 skorobogatov09descent, demarche09obstruction, poonen10insufficiency,
 cdx19comparing} and
\eqn{
\XA^\ddesc=\XA^{\fdesc}=\XA^\desc
}
 \cite{cao20sous}.
See Section  \ref{ptob} for the above notations.
Moreover, for smooth geometrically connected variety,
 Harpaz and Schlank \cite{hs13homotopy} defined the (resp. linearised)
 homotopy obstruction  $\XA^{h}$ (resp. $\XA^{\ZZ h}$ and showed that
 $\XA^{\Br}=\XA^{\ZZ h}$ and $\XA^{\etBr} = \XA^{h}$.

A natural question is, do these relations hold on algebraic stacks ?
Recently, the author and Wu \cite{lh23stackbm} showed that
 $\XA^{\Br} \subseteq \XA^{\conn}$ for any algebraic $k$-stack $X$.
They also showed that $\XA^\ddesc=\XA^\desc$ for a quotient stack
 $X = [Y/F]$  where  $Y$ is a quasi-projective smooth geometrically integral
 $k$-variety, and $F$ is a finite $k$-group
 \cite{wl24stackdd}.

In the papers mentioned above,
 the methods
 to extend relations from varieties to algebraic stacks
 are quite \emph{ad hoc}.
They could not deal with more than one relations at once.

The most suitable way is cohomological descent,
 but is  not  good at dealing with individual degrees of cohomology
 groups (such as $\Br$), and {\adelic} points.

In this paper, we use the advantage of cohomological descent and six
 operation formalism in the infinity categorical   setting
 (Liu and Zheng \cite{lz17enhanced}),
 to develop a systematical way for extending relations to algebraic
 stacks.
The methods in  \cite{lz17enhanced} enhanced the usual derived category
 into closed symmetric monoidal presentable stable $\infty$-categories
 encoding six operations and base chage theorems, {\Kunn} formula,
 {\Poin} duality (and so on), such that they satisfies cohomological
 descent in a homotopy coherent way.
We use this formalism  to construct certain
 $\infty$-categories, carrying data of
 {\adelic} points  and obstructions,
 and show that they
  also satisfies cohomological
 descent.
In the case that the obstruction is defined by  abelian cohomology
 groups (such as $\Br$), the $\infty$-categories above are
 able to completely recover the obstruction sets.
As an application, we extend several relations between obstructions
 from varieties to algebraic stacks, under mild assumptions.

The main results are described as follows.
\thmu{ [{Thm. \ref{thm_tobs}}]
Let $\Sch_{/k}'$ be a full subcategory of the category of $k$-schemes,
 $\Esp_{/k}'$ (resp. $\Chp_{/k}'$) be the corresponding
 full subcategory of the category of algebraic $k$-space
 (resp. sub $2$-category of the $2$-category of algebraic $k$-stacks)
 build from $\Sch_{/k}'$ (resp. $\Esp_{/k}'$)-atlas.
Let $\sE$ be a  family of smooth surjective representable maps of $\Chp_{/k}'$
 satisfying some mild conditions, and
\eqn{
\obs: \Sch_{/k}'\ra \Set, \quad X\mpt \XA^\obs
}
 be a map.
Then to
 every pair $(\obs, \sE)$
 we  may  associate to
 a subset $X(A)^{\tobs_\sE}$ such that  $\tobs_\sE$
 defines an obstruction on ${\Chp_{/k}'}^\sE$ (whose objects
 are build from $\sE$-atlas, see Notation \ref{nota_F1_F2})
 satisfying the following properties.
\enmt{[\upshape (1)]
\item \label{uit_tobs_A}
For any $X\in \Sch_{/k}'$,
 we have $\XA^\obs\subseteq \XA^{\tobs_\sE}$.

\item \label{itu_tobs_func}
We have  ${\tobs_\sE}$ is functorial on $({\Chp_{/k}'}^\sE)_\rep$
 $($Definition \ref{defi_obs_functorial}$)$.

\item \label{itu_tobs_cmp1}
Suppose that we are given two maps
\eqn{
\obs_1\subseteq \obs_2: \Sch_{/k}'\ra \Set
}
Then for every $X\in {\Chp_{/k}'}^\sE$,
 we have $\XA^{\wt{\obs_1}_\sE}\subseteq \XA^{\wt{\obs_2}_\sE}$.

\item \label{itu_tobs_B}
If
 for each $f: Y\ra X$ in $\sE\cap (\Sch_{/S})_1$,
  $f$  induces a surjective map
  $Y(A)^\obs\ra X(A)^\obs$,
 then we have the following properties.
\enmt{[\upshape (i)]
\item \label{itu_tobs_only}
Let  $X\in {\Chp_{/k}'}^\sE$ such that there is a
 a map
 $f: Y\ra  X$
 in $({\Chp_{/k}'}^\sE)_\rep$ with $Y\in \Sch_{/k}'$,
 and  $\obs$ is the only one for  $Y$
 $($i.e., $\YA^\obs\neq\emptyset$ implies
 $Y(k)^\obs\neq\emptyset)$,
 then  $\tobs_\sE$ is
 the only one for  $X$.

In particular, if $\obs$ is the only one for every $X\in \Sch_{/k}'$,
 then  $\tobs_\sE$ is
 the only one for every $X\in {\Chp_{/k}'}^\sE$.

\item \label{itu_tobs_Sch}
  If moreover,  $\obs$ comes from a cohomological functor (Definition
    \ref{defi_B0_fai_coh_fun})
    $\obs: \Chp_{/k}'^\op\ra \Set$,
 then
 for every $X\in {\Chp_{/k}'}^\sE$,
 we have $\XA^{\tobs_\sE}\subseteq \XA^\obs$.

In particular,
 for any $X\in \Sch_{/k}'$,
 we have $\XA^\obs= \XA^{\tobs_\sE}$.
}

\item \label{itu_tobs_E1}
Let $\sE_0\subset\sE$ be another family satisfying conditions
 at the beginning of the theorem.
Then for any $X\in {\Chp_{/k}'}^{\sE_0}$,
 we have $\XA^{\tobs_{\sE_0}}\subseteq \XA^{\tobs_\sE}$.

\item \label{itu_tobs_Fobs_Chp}
Suppose that
 $\obs$ extends  to a map
\eqn{
  \obs: \Chp_{/k}'^\sE\ra \Set, \quad X\mpt X(A)^\obs\subseteq X(A).
}
Let  $X\in {\Chp_{/k}'}^\sE$ such that there is
 a map
 $f: Y\ra  X$
 in $({\Chp_{/k}'}^\sE)_\rep$ with $Y\in \Sch_{/k}'$,
  and   $f$  induces a surjective map
  $Y(A)^\obs\ra X(A)^\obs$,
Then   we have $\XA^\obs \subseteq\XA^{\tobs_\sE}$.

In particular, if for
  each $f: Y\ra X$ in $({\Chp_{/k}'}^\sE)_\rep$ with $Y \in \Sch_{/k}'$,
  $f$  induces a surjective map
  $Y(A)^\obs\ra X(A)^\obs$,
 then
 for every $X\in {\Chp_{/k}'}^\sE$,
 we have $\XA^\obs \subseteq\XA^{\tobs_\sE}$.
If moreover $\obs$ comes from a cohomological functor
    $\obs: \Chp_{/k}'^\op\ra \Set$,
 then
 for every $X\in {\Chp_{/k}'}^\sE$,
 we have $\XA^\obs= \XA^{\tobs_\sE}$.

\item \label{itu_tobs_ttobs}
For every $X\in {\Chp_{/k}'}^\sE$,
 we have $\XA^{\wt{(\tobs_\sE)}_\sE} = \XA^{\tobs_\sE}$.
}}

The paper is organised as follows.
Section \ref{ptob} gives a functorial way to extend
 definition of classical obstructions (mostly cohomological ones)
 from varieties to algebraic stacks.
Then we use six operation formalism
 and enhanced derived categories to construct  $\infty$-categories
 mentioned before encoding obstructions in Section \ref{six}.
While Section \ref{cohdes} is devoted to prove cohomological descent
 properties of the constructed  $\infty$-categories.
In Section \ref{appl}, we apply these formalism to
 construct new obstructions, which
 simultaneously  ``extend" and ``approximate" existing obstructions
 from varieties to algebraic stacks,
At last, we give some concrete examples in Section \ref{conc} to
 illustrate what these new obstruction like.

\section{Points and obstructions on algebraic stacks} \label{ptob}
Let us briefly introduce basic notions of rational and {\adelic}
 points, obstructions  to local-global principle on algebraic stacks.

We omit set theoretic issues since we can carefully avoid it by using
 {\Grot} universes.
Denote by $\Set$ the category of sets, $\Sch$ the category of
  quasi-separated
 schemes,
 $\Esp$ the category of algebraic spaces and $\Chp$ the $(2, 1)$-category of
 algebraic stacks.
Here we talk about algebraic spaces and stacks
 in the sense of \cite[025Y,026O]{sp}, which is based on big fppf
 topology.

For a commutative ring $R$, we also write $R$ for $\Spec R$.
Fix a base scheme $S\in \Sch$ and let $X, T\in \Ob(\Chp_{/S})$ be two
  algebraic stacks over $S$.
The \emph{$T$-points} of $X$ is the set of isomorphism classes
 of objects in the groupoid $\Hom_{\Chp_{/S}} (T, X)$, denoted by $X(T)$.
Let $k$ be a global field with {\adele} ring $\bfA_k$,
 the natural inclusion $k\subset \bfA_k$ induces $q:  \Spec \bfA_k \ra \Spec k$,
 making $\Spec \bfA_k$ an object of $\Chp_{/k}$.
Let $X\in\Ob(\Chp_{/k})$.
Then we call  $\XA$ (resp. $\Xk$)  the \emph{{\adelic} points} (resp.
  \emph{rational points}) of $X$.

\rk{
If $X$ is represented by a scheme, this coincides with the classical
 definition.
But in general, for $X$ being an algebraic stack,
 the map  $\Xk\ra \XA$ induced by $q$
 is not necessary injective.
For example, let $G$ be an affine $k$-group.
Let $*$  be the neutral element in  the pointed set
 $BG(\bfA_k)) = H_\fppf^1(\bfA_k, G)$.
Then its preimage $\ker(BG(k)\ra BG(\bfA_k))$
 defined by the following {\Cart} diagram
\eqn{\xymatrix{
&\ker(BG(k)\ra BG(\bfA_k))\ar@{=}[ddl] \ar[r] \ar@{^(->}[d]
 &\{*\}\ar@{^(->}[d] \\
 & BG(k)\ar[r]\ar@{=}[d] & BG(\bfA_k) \ar@{=}[d]\\
\Sha^1(G/k)\ar[r]& H^1(k, G)\ar[r] &\check H_\fppf^1(\bfA_k, G)
}}
  is   the {\TS}   group $\Sha^1(G/k)$, which
 is not necessarily trivial.

Anyway, by abuse notation, we also regard $\Xk$ as the image of
 $\Xk\ra\XA$ and write $\Xk\subseteq \XA$.
}

\defi{  [{\cite[2.17]{lv2desc}}] \label{defi_stable}
Let $\sC$ be a $(2, 1)$-category, and $\sD$ be an  ordinary
 category.
Let $F: \sC\ra\sD$ be a functor
 from the underlying ordinary category of $\sC$ to $\sD$.
We say that $F$ is \emph{stable} if
 $F$ is a strict $2$-functor from $\sC$ to $\sD$,
 i.e., for any $1$-morphisms $f$ and $g$ in  $\sC$
 that is $2$-isomorphic, we have $F(f)=F(g)$.
}
\rk{
There is a notion of \emph{stable functor} in the literature
 which means ``having a left adjoint on each slice".
In this paper, we  keep the
  nomenclature in Definition \ref{defi_stable}.
}

Let $X\in \Chp_{/k}$.
For a stable functor $F: (\Chp_{/k})^\op\ra \Set$ and $\alpha\in F(X)$,
 as the classical case in Poonen \cite[8.1.1]{poonen17rational},
 one may also define
 the \emph{obstruction given by $\alpha$} to be  the subset
 $\XA^\alpha$ of $\XA$ whose elements are characterized by
\eqn{
\XA^\alpha=\{x\in \XA\mid \alpha(x)\in \im F(q)\},
}
 (which is well-defined by stability of $F$),
 and the  \emph{$F$-set} (or \emph{$F$-obstruction})
 to be the subset $\XA^F$ whose elements are
 characterized by
\eq{ \label{eq_XAF}
\XA^F=\bigcap_{\alpha\in F(\cX)} \XA^\alpha =
 \{x\in \XA\mid \im F(x)\subseteq \im F(q)\}.
}
Then we have the inclusion
\eqn{
 \Xk\subseteq \XA^F\subseteq  \XA^\alpha\subseteq \XA.
}

\rk{
There is a more functorial way to define $\XA^F$.
Namely, one sees that the map $\Xk\ra \XA$ obviously
 factorizes as
\eq{ \label{eq_obs}
\Xk\ra \Map(F(X), F(k))\tm_{F(q), \Map(F(X), F(\bfA_k)), F}
 \XA\xra{q^*} \XA,
}
 and we define  $\XA^F = \im(q^*)$.
Also let $A\in F(X)$, if we replace $F(X)$ in \eqref{eq_obs}
 by $\{A\}$,
 then we define   $\XA^A$ to be the resulting set.
It is easy to verify that these definitions  coincide
 with those by \eqref{eq_XAF}.
}

Now let $F = H_\et^n(-, G)$ for some commutative $k$-group $G$.
Then $F: (\Chp_{/k})^\op \ra \Set$
 is stable (c.f. \cite[2.29]{lv2desc})
 and the
 resulting $F$-set  $\XA^{H_\et^n(-, A)}$ coincides with classical definition
 if  $X$ is a $k$-variety (by a $k$-\emph{variety} $X$ we mean an object of
 $\Sch_{/k}$ that is separated of finite type).

In particular, $F = \Br = H_\et^2(-, \bfG_m)$ gives
 the \emph{{\BM} obstruction} of $\XA^{\Br}$.

Next let $G$ be a $k$-group and  $T\in\Ob(\Chp_{/k})$, write $G_T=G\tm_k T$.
Then one considers fppf
 $G_T$-torsors
 over $T$
 (see, for example, \cite[Sec. 2.3]{lh23stackbm}).
There is an isomorphism of pointed sets (c.f.
 \cite[III.3.6.5 (5), IV.3.4.2 (i)]{giraud71cohnonab})
\eqn{
\check H_\fppf^1(X, G)\xra{\sim} \{\text{$G$-torsors over $X$}\}/\cong.
}
It turns out that
 $F=\check H_\fppf^1(-, G): (\Chp_{/k})^\op\ra \Set$ is also a stable functor
 (c.f. \cite[3.14]{lv2desc})
 and the
 resulting $F$-set  $\XA^{\check H_\fppf^1(-, G)}$ coincides with classical
 definition
 if
 $X$ is a $k$-variety. See \cite[3.17 (i)]{lv2desc}.
Define  the \emph{descent obstruction} of $X$ to be
\eqn{
\XA^\desc=\bigcap_{\text{linear $k$-group $G$}} \XA^{\check H_\fppf^1(-, G)}.
}
Similarly, define
\gan{
\XA^{\PGL}=
 \bigcap_{n\ge 1} \XA^{\check H_\fppf^1(-, \PGL_n)}, \\
\XA^\conn=\bigcap_{\text{connected linear $k$-group $G$}}
 \XA^{\check H_\fppf^1(-, G)},
}
 and the \emph{second descent obstruction}
   \cite[4.2]{lv2desc}
\eqn{
\XA^{\sdesc}=
 \bigcap_{\text{commutative $k$-group $G$}}
 \XA^{H_\et^2(-, G)}.
}

Let $f: Y\xra{G} X$ be a $G$-torsor over $X$
 (which is also in $\Chp_{/k}$ by \cite[Lem. 2.10]{lh23stackbm})
 and denoted by $f$
 the element of
 $\check H_\fppf^1(X, G)$ corresponding to the class of
 $f: Y\xra{G} X$.
\prop{
We have
\eqn{
\XA^f=\bigcup_{\s\in H^1(k, G)}
 f^\s(Y^\s(\bfA_k)),
}
 where $f^\s: Y^\s\xra{F^\s} X$  is the twist of
  $f$ $($see \cite[3.30]{lv2desc}$)$.
}
\pf{
For $X$ being a  $k$-variety,  this is well-known.
For algebraic stacks, this comes from general descent by torsors
 \cite[3.20]{lv2desc}.
}
From this we may define  the \emph{{\etale} Brauer  obstruction} of $X$, to be
\eqn{
\XA^\etBr=\bigcap_{f:Y\xra{F} Y \text{ torsor under
 finite $k$-group $F$}} \bigcup_{\s\in H^1(k, F)}
 f^\s(Y^\s(\bfA_k)^{\Br}),
}
 the \emph{finite  descent  obstruction}
\eqn{
\XA^\fdesc=\bigcap_{f:Y\xra{F} Y \text{ torsor under
 finite $k$-group $F$}} \bigcup_{\s\in H^1(k, F)}
 f^\s(Y^\s(\bfA_k)^{\desc}),
}
 and the \emph{iterated descent  obstruction} of $X$, to be
\eqn{
\XA^\ddesc=\bigcap_{f:Y\xra{G} Y \text{ torsor under
 linear $k$-group $G$}} \bigcup_{\s\in H^1(k, G)}
 f^\s(Y^\s(\bfA_k)^{\desc}).
}

All these  definitions of obstructions coincide with  the
 classical ones in the case that $X$ is a $k$-variety.

\section{Six operation formalism and categories encoding obstructions}
 \label{six}

Let $X$ be an topos and  $\La$ be a sheaf of rings in $X$.
Recall that $\sD(X, \La)$ is the
 closed symmetric monoidal presentable stable $\infty$-category constructed in
 \cite[Sec. 2.2]{lz17enhanced}.
It is an enhanced version of the classical derived categories.
More specifically,
 let    $\PTop$ be the $(2,1)$-category of topos and $\Rind$ the category
 of ringed diagrams (\cite[Def. 2.2.5]{lz17enhanced}).
We have a functor
\eq{ \label{eq_D_topos}
\sD: \N(\PTop)^\op\tm \N(\Rind)^\op \ra \Cat_\infty
}
 which sends a map $f: Y \ra X$ in $\PTop$ and $\lam\in\Rind$
 to $f^*: \sD(X, \lam)\ra \sD(Y, \lam)$ (c.f.
  \cite[Notation 2.2.6]{lz17enhanced}).

Now fix a final object $e$
 and let $u_0: U_0\ra e$ be a covering in $X$ with {\Cech} nerve
 (see the definition after \cite[Prop. 6.1.2.11]{HTT})
 $u_\bu: U_\bu\ra e$.
Put $\La_n= \La\tm U_n$.
For any $\alpha: [m]\ra [n]$ in $\bDel_+$,  we have
 $u_\alpha^*:  \sD(X_{/U_m}, \La_m) \ra \sD(X_{/U_n}, \La_n)$.
Then  we have an augmented
 cosimpicial $\infty$-category $\sD(X_{/U_\bu}, \La_\bu)$.
It has the following descent property.
\lemu{[{\cite[Lem. 3.1.3]{lz17enhanced}}] \label{lemm_D_lim}
The natural map
\eqn{
\sD(X, \La)\ra \varprojlim_{n\in\bDel}\sD(X_{/U_n}, \La_n)
}
 is an equivalence of $\infty$-categories.
}

\lemm{ \label{lemm_D_rad}
For any $\alpha: [m]\ra [n]$ in $\bDel_+$,
 the diagram
\eqn{
\xymatrix{
\sD(X_{/U_m}, \La_m)\ar[r]^-{u_{d_0^{m+1}}^*}\ar[d]^-
 {u_\alpha^*} &\sD(X_{/U_{m+1}}, \La_{m+1})\ar[d]^-{u_{\alpha'}^*} \\
\sD(X_{/U_n}, \La_n)\ar[r]^-{u_{d_0^{n+1}}^*}\ &\sD(X_{/U_{n+1}}, \La_{n+1}).
}}
is right adjointable $($see \cite[Def. 1.4.1]{lz17enhanced}$)$.
}
\pf{
Apply \cite[Lem. 2.2.9]{lz17enhanced} to
the transpose of the square.
}

Now let $A$ be another topos.
Fix $\La$ to be a ring, representing a constant sheaf of rings in
 any topos.
\defi{ \label{defi_H_topos}
Let $\sH(X, \La)=\Fun(\sD(A, \La), \sD(X, \La))$.
We also have a functor
\eq{ \label{eq_H_topos}
\sH:  \N(\PTop)^\op \ra \Cat_\infty
}
 which maps $f: Y\ra X$ to $f^*: \Fun(\sD(A, \La), \sD(X, \La))\ra
  \Fun(\sD(A, \La), \sD(Y, \La))$.
}
By \cite[Prop. 5.5.3.6]{HTT} 
 $\sH(X, \La)$ is a
 presentable $\infty$-category.

\prop{ \label{prop_H_lim}
The natural map
\eqn{
\sH(X, \La)\ra \varprojlim_{n\in\bDel}\sH(X_{/U_n}, \La)
}
 is an equivalence of $\infty$-categories.
}
\pf{
For any $\infty$-category $\sC$, the functor
 $\Fun(\sC, -): \Cat_\infty \ra \Cat_\infty$ has a left adjoint
 $-\tm \sC: \Cat_\infty \ra \Cat_\infty$.
Thus  $\Fun(\sC, -)$ preserves all small
 limits (\cite[Prop. 5.2.3.5]{HTT}).
The result follows from Lemma \ref{lemm_D_lim}.
}

\nota{ \label{nota_F1_F2}
Let  $\sF_1$ and $\sF_2$ be two family of maps of $\Chp_{/S}$
 that are stable under composition and pullback and both contain
 every degenerate edge.
Then we denote by  $(\Esp_{/S}^{\sF_1})_{\sF_2}$ the subcategory of
 $\Esp_{/S}$ spanned by objects that admits
 $\sF_1\cap (\Esp_{/S})_1$-atlas,
 that is, those $X\in \Esp_{/S}$ admits a {\etale} surjective
 map $Y\ra X$ in  $\sF_1\cap (\Esp_{/S})_1$ such that $Y\in \Sch_{/S}$,
 and by morphisms in $\sF_2\cap (\Esp_{/S})_1$.

Similarly, denote by  $(\Chp_{/S}^{\sF_1})_{\sF_2}$ the sub $2$-category of
 $\Chp_{/S}$ spanned by objects that admits $\sF_1$-atlas from
 $\Esp_{/S}^{\sF_1}$,
 that is, those $X\in \Chp_{/S}$ admits a smooth surjective
 map $Y\ra X$ in  $\sF_1$ such that $Y\in \Esp_{/S}^{\sF_1}$,
 and by $1$-morphisms in $\sF_2$.

Let $\Box$ be the family of maps $Y\ra X$ such that $X$ is
  $\Box$-coprime (that is, there exists a map $X\ra \Spec \ZZ[\Box^{-1}]$).
Then clearly $Y$ is also $\Box$-coprime.
Let
 $\rep$ be the family of maps that are representable,
 $(\sr)$ be the family of maps that are smooth representable,
 and $(\ssr)$ be the family of maps that are smooth surjective representable.
}

Recall in \cite[Def. 5.5.3.1]{HTT} that
 $\PrL$ is the subcategory of $\Cat_\infty$ spanned by
 presentable $\infty$-categories and functors that preserve small colimits.
By \cite[Cor. 5.5.2.9]{HTT}, a functor between two presentable
 $\infty$-categories has a right adjoint
 if and only if
 it preserves small colimits,
 i.e. it is an edge of  $\PrL$.
We will freely use this fact in the rest of this paper.

Now we work in $\Chp_{/S}$.
Recall in \cite[(5.7)]{lz17enhanced},
 we have a functor
\eq{ \label{eq_D}
\sD: \N(\Chp_{/S})^\op \ra \PrL
}
 (after composing  the forgetful functor $\N(\Chp_{/S})\ra \N(\Chp)$)
 such that,
 it coincides on $\N(\Chp_{/S})_\sr^\op$
 with the functor \eqref{eq_D_topos} compositing
 $X\mpt X_\liset$
 (\cite[Cor. 5.3.8]{lz17enhanced}), and
 for $X\in \Chp_{/S}$,
 we have $\sD(X, \La)$
 whose homotopy category
 $h\sD(X, \La)=D_\cart(X_\liset, \La)$ (\cite[5.3.8]{lz17enhanced}).
For convince, we write $\sD(X)$ for $\sD(X, \La)$,
 $D(X)$ for $D_\cart(X_\liset, \La)$.

For any map $f: Y\ra X$ in $\Chp_{/S}$ we have operations
\eq{ \label{eq_D_adj_*}
f^*: \sD(X)\rightleftarrows \sD(Y): f_*,
}
 and for $f$ is locally of finite type, $\La\in\Ring_{\Box\tu{-tor}}$
  and $X$ is $\Box$-coprime,
\eq{ \label{eq_D_adj_!}
f_!: \sD(Y)\rightleftarrows \sD(X): f^!,
}
 where $\Box$ is a nonempty set  of rational primes,
 $\La\in\Ring_{\Box\tu{-tor}}$ is the category of $\Box$-torsion rings.
See \cite[Sec. 6.2]{lz17enhanced}.

If $f$ is smooth, by {\Poin} duality \cite[Thm. 6.2.9]{lz17enhanced},
 $f^*\langle\dim f\rangle\xra{\sim} f^!$.
In particular if  $f$ is {\etale},
 $f^*\xra{\sim} f^!$, i.e.,  $f_!\dashv f^*$.

Fix $A\in\Chp_{/S}$ with structure map  $q: A\ra S$,
 we can make a similar definition as in Definition \ref{defi_H_topos},
 namely,
\eqn{
\sH(X)=\sH(X, \La)=\Fun(\sD(A, \La), \sD(X, \La)).
}
It is presentable.
By \eqref{eq_D}, this produces   a functor
\eqn{ \label{eq_H}
\sH: \N(\Chp_{/S})^\op \ra \PrL,
}
 which coincides on $\N(\Chp_{/S})_\sr^\op$
 with the functor \eqref{eq_H_topos} compositing
 $X\mpt X_\liset$. See comments after \eqref{eq_D}.

Suppose that $f: X_0\ra X_{-1}$ is smooth surjective map in
 $\Chp_{/S}$, and
 $X_\bu: \N(\bDel_+)^\op \ra \N(\Chp_{/S})$ its {\Cech} nerve.
Then
 by \cite[Thm. 6.2.13 (1)]{lz17enhanced},
 the natural map
\eqn{
\sD(X_{-1})\ra \varprojlim_{n\in\bDel}\sD(X_n)
}
 is an equivalence of $\infty$-categories.
Similar to the proof  of
 Proposition \ref{prop_H_lim}, we obtain
  an induced equivalence
\eqn{ 
\sH(X_{-1})\ra \varprojlim_{n\in\bDel}\sH(X_n).
}

\defi{
Suppose that the structure map of $X$ is $p: X\ra S$.
Define   $\infty$-categories  $\sH_0(X)$ and $\sH_1(X)$ by the {\Cart} diagram
\eq{ \label{eq_H0_H1}
\xymatrix{
\sH_1(X)\ar[r]^-{V}\ar[d]^-{Q_*}
   &\sH_0(X)\ar[r]^-{U}\ar[d]^-{P_*} &\sH(X)\ar[d]^-{p_*} \\
 \sH(S)_{q_*\ds q_*}\ar[r]^-{v}  &\sH(S)_{q_*/}\ar[r]^-{u}  &\sH(S)
}}
 where  $u$ and $v$ is the natural forgetful functor,
 and  we write $\sH(S)_{q_*\ds q_*}$ for
 $(\sH(S)_{q_*/})_{/(q_*\xra{\id}q_*)}$.
}

\rk{ \label{rk_H0_H1_struc}
\enmt{[\upshape (1)]
\item \label{it_H0_u}
We already know that $\sH(X)$ and $\sH(S)$  are presentable.
By \cite[Prop. 5.5.3.11]{HTT}, so is  $\sH(S)_{q_*/}$.
In particular, they all admits small limits \cite[Prop. 5.5.2.4]{HTT}.
By dual version of \cite[Prop. 1.2.13.8 (1)]{HTT},
 such limits are preserved by $u$.

\item \label{it_H0_H1_struc}
Since $u$ (resp. $v$) is a left (resp. right) fibration by
 (resp. dual version of)
 \cite[Cor. 2.1.2.2]{HTT}.
Then $u$ and $v$ are all categorical fibration,
 and hence \eqref{eq_H0_H1} is also homotopy {\Cart}.
With loss of generality, one may view
 $\sH_0(X)$ (resp. $\sH_1(X)$) as the $\infty$-category
 whose vertices consist of
 pairs  $(x_*, \lam_*)$
 (resp. triples $(x_*, \lam_*, \mu_*)$),
 where $x_* \in \sH(X)_0$,
 $\lam_*: q_*\ra p_*x_*$ in $\sH(S)_1$ (resp. and
 $\mu_*: p_*x_*\ra q_*$ in $\sH(S)_1$ satisfying
  $\mu_*\lam_*\xra{\sim}\id_{q_*}$).
An edge
 $({x_1}_*, {\lam_1}_*)\ra
 ({x_2}_*, {\lam_2}_*)$ of $\sH_0(X)$
 (resp.
 $({x_1}_*, {\lam_1}_*, {\mu_1}_*)\ra
 ({x_2}_*, {\lam_2}_*, {\mu_2}_*)$ of $\sH_1(X)$)
 is a
 pair $(\xi_*, \s)$ (resp.
 a triple $(\xi_*, \s_1, \s_2)$) where
 $\xi_*$ is an
 edge $\xi_*: {x_1}_*\ra   {x_2}_*$ in $\sH(X)$ and
  $\s: \Delta^2\ra \sD(S)$ is a map (resp.
  $\s_1, \s_2: \Delta^2\ra \sD(S)$ are two maps) representing
  a homotopy $(p_*\xi_*)\circ{\lam_1}_* \sim {\lam_2}_*$
  (resp. two homotopies
 $(p_*\xi_*)\circ{\lam_1}_* \sim {\lam_2}_*$ and
 ${\mu_2}_*\circ\xi_* \sim p_*{\mu_1}_*$).

The functor $U$ (resp. $V$) maps
   $(x_*, \lam_*)$ (resp. $(x_*, \lam_*, \mu_*)$)
    to
  $x_*$ (resp.   $(x_*, \lam_*)$),
  and the functor $P_*$ (resp. $Q_*$) mpas
   $(x_*, \lam_*)$ (resp. $(x_*, \lam_*, \mu_*)$)
   to
   $\lam_*: q_*\ra p_*x_*$
   (resp. $q_*\xra{\lam*} p_*x_*\xra{\mu_*} q_*$).

\item \label{it_H0_V}
Since  $v$ is a ritht fibration.
By   \cite[Cor. 2.4.2.3 (2)]{HTT}, so is $V$.
In particular $V$ is a {\Cart} fibration.
Thus  given $({x_2}_*, {\lam_2}_*, {\mu_2}_*)\in \sH_1(X)_0$,
 an edge
 $({x_1}_*, {\lam_1}_*)\ra
 ({x_2}_*, {\lam_2}_*) = V({x_2}_*, {\lam_2}_*, {\mu_2}_*)$ of $\sH_0(X)$
 can be lifted to an edge
 $({x_1}_*, {\lam_1}_*, {\mu_1}_*)\ra
 ({x_2}_*, {\lam_2}_*, {\mu_2}_*)$ of $\sH_1(X)$
  which
  maps to
 $({x_1}_*, {\lam_1}_*)\ra
 ({x_2}_*, {\lam_2}_*)$
  under $V$.
}}

For $X\in \Chp_{/S}$,
 we  denote by $p_X$ the structure map $X\ra S$  of $X$.
If there is
 no confusing,  we simply write $p=p_X$.

\lemm{ \label{lemm_H0_adj}
Suppose that $f: Y\ra X$ is map in $\Chp_{/S}$,
Then the adjoint pair
\eqn{
f^*: \sH(X)\rightleftarrows \sH(Y): f_*
}
 indeced by \eqref{eq_D_adj_*}
 induces an adjoint pair
\eqn{
f^*: \sH_0(X)\rightleftarrows \sH_0(Y): f_*.
}}
\pf{
Since we have an the natural equivalence
 $\ep: {p_Y}_*\xra{\sim}{p_X}_*f_*$
 deduced from the functoriality of $\sD$,
 the fact that $f_*$ induces a functor  $\sH_0(Y)\ra \sH_0(X)$
 is clear, sending $(y_*, \lam_*)$ to $(f_*y_*, (y_*\ep)\lam_*)$.
To  show that  $f^*$ induce a functor $\sH_0(X)\ra \sH_0(Y)$
 we use the unit map $a: \id\ra f_*f^*$,
 namely, we construct the disired functor
 by sending
 $(x_*, \lam_*)$ to
 $(f^*x_*, (\ep^{-1}f^*x_*) ({p_X}_*ax_*) \lam_*)$.
Then the checks of the two adjunctions  are straightforward.
The proof is complete.
}
\rk{ \label{rk_H0_pb_exist}
The functor
 $f^*: \sH_0(X)\ra \sH_0(Y)$
 constructed in Lemma \ref{lemm_H0_adj} can also be
 described as follows.
Let ${P_X}_*^a: \sH_0(X)\ra \sH(S)_{p_*/}$ be the functor sending
  $(x_*, \lam_*)$ to $({p_X}_*ax_*)\lam_*: q_*\ra {p_X}_*f_*f^*$.
Then we have the homotopy commutative diagram with  {\Cart} square
\eqn{
\xymatrix{
  \sH_0(X)\ar@/^/[drr]^-{f^*\circ U_X} \ar@{.>}[dr]|-{f_0^*}
   \ar@/_/[ddr]_-{{P_X}_*^a} & & \\
   &\sH_0(Y)\ar[r]_-{U_Y}\ar[d]^-{{P_Y}_*} &\sH(Y)\ar[d]^-{{p_Y}_*} \\
   &\sH(S)_{q_*/}\ar[r]^-{u}  &\sH(S).
}}
It follows that there is a unique functor $f_0^*$ rendering the diagram
 commutitive.
Then one checks that $f_0^*\xra{\sim} f^*$ is the one constructed above.
}

\lemm{ \label{lemm_lim_pullback}
Let
\eqn{
\xymatrix{
\sX'\ar[r]^-{q'}\ar[d]^-{p'} &\sX\ar[d]^-{p} \\
 \sY'\ar[r]^-{q}  & \sY
}}
 be a diagram of $\infty$-categories which is homotopy {\Cart}
 $($with respect to the Joyal model structure$)$ and let $K$
 be a simplicial set.
Suppose that $\sX$ and $\sY'$ admit limits for all diagrams indexed
 by $K$ and that $p$ and $q$ preserve
 limits of diagrams indexed by $K$.
Then
\enmt{[\upshape (1)]
\item \label{it_pullback_is_lim}
 A diagram $\ol f: K^\tril\ra \sX'$ is a limit of $f = \ol f|K$
  if and only if $p'\ol f$ and $q'\ol f$ are limit diagrams.
 In particular, $p'$ and $q'$ preserve limits
 indexed by $K$.
\item \label{it_pullback_adm_lim}
Every diagram $f: K \ra \sX'$ has a limit in $\sX'$.
}}
\pf{
We have the equivalence $R: \Cat_\infty\ra \Cat_\infty$
  associating to every $\infty$-category its opposite
 (\cite[Rmk. 2.4.2.7]{HA}).
Then apply $R$ to \cite[Lem. 5.4.5.5]{HTT}. 
}

\lemm{ \label{lemm_H0_U_lim}
Let $X\in \Chp_{/S}$.
\enmt{[\upshape (1)]
\item \label{it_H0_adm_lim}
The category $\sH_0(X)$ admits small limits.
\item \label{it_U}
Let  $K$ be a   small simplicial set and $a: K\ra \sH_0(X)$ be a diagram.
Then  $K^\tril\xra{a} \sH_0(X)$ is a limit diagram
 if and only if the composition
 $K^\tril\xra{a} \sH_0(X)\xra{U} \sH(X)$
 is a limit diagram, where $U$ is the functor defined in \eqref{eq_H0_H1}.
}}
\pf{
By Remark  \ref{rk_H0_H1_struc} \eqref{it_H0_u},
 we  know that $\sH(X)$, $\sH(S)$
  and $\sH(S)_{q_*/}$ admits small limits,
  and such limits are preserved by $u$.
Since $p_* \in \PrR$, it also
 preserves small limits \cite[Prop. 5.2.3.5]{HTT}.
Then \eqref{it_H0_adm_lim}  follows from Lemma
 \ref{lemm_lim_pullback} \eqref{it_pullback_adm_lim}.

For \eqref{it_U},
 I claim that    $K^\tril\xra{a} \sH_0(X)\xra{U}\sH(X)$ is
 a limit diagram implies that $K^\tril\xra{a} \sH_0(X)\xra{P_*}
 \sH(S)_{q_*/}$ is a limit diagram.
Indeed,
 by  the dual version of \cite[Prop. 1.2.13.8 (2)]{HTT},
  we only need to show that
 $K^\tril\xra{a}\sH_0(X)\xra{U}  \sH(X)\xra{p_*} \sH(S)$
 is a limit diagram, which  is clear since
 $p_*$ preserves small limits.
Then  we use  Lemma
 \ref{lemm_lim_pullback} \eqref{it_pullback_is_lim}
 to obtain the desired result.
}

\constr{ \label{constr_coh}
Let $\ol x_* = (x_*, \lam_*, \mu_*)
 \in \sH_1(X)$ (see Remark \ref{rk_H0_H1_struc} \eqref{it_H0_H1_struc}).
Suppose that $x_*$ has a left adjoint
 $x^*: \sD(X)\ra \sD(A)$.
Then by definition of $\sH_1(X)$,  we have edges
\eqn{
q^*\xra{\mu^*} x^*p^*\xra{\lam^*} q^*
}
 of $\Fun(\sD(X), \sD(A))$.
Fixing tow vertices $K_1, K_2\in \sD(S)$, we also view them
 as  objects in the homotopy
 category $h\sD(S) = D(S)$.
For  an edge of the form $\alpha: p^*K_1\ra p^*K_2$ in $\sD(X)$,
 the composition
\eqn{
q^*K_1\xra{\mu^*K_1} x^*p^*K_1\xra{x^*\alpha} x^*p^*K_2\xra{\lam^*K_2} q^*K_2
}
 is an edge of $\sD(A)$.
Passing to the homotopy categories $D(X)$ and $D(A)$,
 we obtain a well-defined homomorphism
\eq{ \label{eq_x^*_D}
\Hom_{D(X)}(p^*K_1, p^*K_2) \ra \Hom_{D(A)}(q^*K_1, q^*K_2).
}
We simply denote this corresponding
 homomorphism by $h\ol x^*$.
}

\lemm{ \label{lemm_hx^*_id}
Let
\eqn{
\ol\xi_*=(\xi_*, \s): \ol{x_1}_* = ({x_1}_*, {\lam_1}_*, {\mu_1}_*)\ra
 \ol{x_2}_* = ({x_2}_*, {\lam_2}_*, {\mu_2}_*)
}
 be an edge of $\sH_1(X)$
 $($see Remark \ref{rk_H0_H1_struc} \eqref{it_H0_H1_struc}$)$
 and suppose that
 both ${x_1}_*$ and ${x_2}_*$ have left adjoints.
Then we have $h\ol{x_1}^* = h\ol{x_2}^*$.
}
\pf{
By assumption we have an edge $\xi^*: x_2^*\ra x_1^*$ of
 $\Fun(\sD(X), \sD(A))$.
By definition of $\sH_0(X)$, we have a map
 $\Delta^3\tm \Delta^1\ra \sD(A)$
 which we depict as follows
\eqn{
\xymatrix{
q^*K_1\ar[r]^-{\mu_2^*K_1}\ar@{=}[d] &x_2^*p^*K_1\ar[r]^-{x_2^*\alpha}
 \ar[d]^-{\xi^*p^*K_1}
 &x_2^*p^*K_2\ar[r]^-{\lam_2^*K_2}
 \ar[d]^-{\xi^*p^*K_2}
 &q^*K_2 \ar@{=}[d] \\
q^*K_1\ar[r]^-{\mu_1^*K_1} &x_1^*p^*K_1\ar[r]^-{x_1^*\alpha}
 &x_1^*p^*K_2\ar[r]^-{\lam_1^*K_2} &q^*K_2,
}}
 whose upper (resp. lower) composition is $h\ol{x_2}^*$ (resp. $h\ol{x_1}^*$).
Then taking homotopy and the result follows.
}

\rk{ \label{rk_coh_id}
In particular,  for $K\in \sD(S)$ and $i\in\ZZ$,
 let $K_1 = \La$ and $K_2 = K[i]$.
The homomorphism \eqref{eq_x^*_D} becomes
 $H_\et^i(X, K)\ra H_\et^i(A, K)$, and we denote the corresponding
 homomorphism $h\ol x^*$ by  $x^*$ if no confusing.
Lemma \ref{lemm_hx^*_id} tells us that for any  edges
  $\ol \xi_* = (\xi_*, \s_1, \s_2):
   \ol{x_1}_*\ra \ol{x_2}_*$ of $\sH_1(X)$,
  we obtain two same homomorphisms
 $x_1^*=x_2^*: H_\et^i(X, K)\ra H_\et^i(A, K)$.
This fact is crucial in the application part  of this paper
 (see Lemma \ref{lemm_B_rec}).
}

\defi{ \label{defi_A}
Let $X, A\in  \Chp_{/S}$ as before.
Let
 $X(A)^\obs\subseteq X(A)$ be a subset.
In particular, we write $\obs = r$ (resp. $\obs=a$) to indicate
 $X(A)^\obs = X(S)$ (resp. $X(A)^\obs = X(A)$)
\enmt{[\upshape (a)]
\item
We define $\sA_1^\obs(X)\subseteq \sH_1(X)$ to be the full
 subcategory spanned by all $\ol x_* = (x_*, \ep, \ep^{-1})
 \in\sH_1(X)_0$,
 where $x\in X(A)^\obs$, and $\ep: q_*\xra{\sim}x_*p_*$ is
 the natural equivalence deduced from the functoriality of $\sD$.
Define  $\sB_1^\obs(X)\subseteq \sH_1(X)$ to be the full
 subcategory spanned by all $\ol x_* \in \sH_1(X)_0$, such that
 there is an edge $\ol x_*\ra \ol{x_0}_*$ of $\sH_1(X)$ where
 $\ol{x_0}_*\in \sA_1^\obs(X)_0$.
\item
Let $\sA_0^\obs(X)$ (resp. $\sB_0^\obs(X)$) $\subseteq \sH_0(X)$
 be the full subcategory spanned by the  essential
 image  of $\sA_1^\obs(X)$ (resp. $\sB_1^\obs(X)$) under the
 forgetful functor
 $V: \sH_1(X) \ra  \sH_0(X)$.
}}

\rk{ \label{rk_A_struc}
By Remark \ref{rk_H0_H1_struc}
 \eqref{it_H0_V}, we have the following facts.
\enmt{[\upshape (1)]
\item \label{it_0_1_ess}
In the definitions of $\sA_0^\obs(X)$ and  $\sB_0^\obs(X)$,
 we may replace essential image by just image.
\item \label{it_B0_lift}
Let  $\ol x_*\ra \ol{x_0}_*$ be an edge of $\sH_0(X)$ where
 $\ol{x_0}_*\in \sB_0^\obs(X)_0$.
Then we have $\ol x_*\in \sB_0^\obs(X)_0$.
}}

\defi{ \label{defi_F}
Suppose that for each $X\in \Chp_{/S}$, we are given a subset
 $X(A)^\obs\subseteq X(A)$.
Then we use $\sFo$ to denote any family of maps of $\Chp_{/S}$  such that
\enmt{[\upshape (a)]
\item \label{it_F_scbc}
The family $\sFo$ is stable under composition and pullback.
\item \label{it_F_cond}
For each $f: Y\ra X$ in $\sFo$,
  $f$ is representable and  induces a surjective map
  $Y(A)^\obs\ra X(A)^\obs$.
}}

\lemm{ \label{lemm_B0_adj}
Suppose that  $f: Y\ra X$ is a map such that $f$ inducing a surjective
 map  $Y(A)^\obs\ra X(A)^\obs$.
Then the adjoint pair
\eqn{
f^*: \sH_0(X)\rightleftarrows \sH_0(Y): f_*
}
 in Lemma \ref{lemm_H0_adj}
 induces an adjoint pair
\eq{ \label{eq_B_adj}
f^*: \sB_0^\obs(X)\rightleftarrows \sB_0^\obs(Y): f_*.
}}
\pf{
Let $y\in Y(A)^\obs$.
Since
 $f: Y\ra X$ induce a
 map  $Y(A)^\obs\ra X(A)^\obs$.
 we know that $fy\in X(A)^\obs$.
Also we have the natural equivalence $(fy)_*\xra{\sim}
 f_* y_*$.
Thus we readily see that
 $f_*: \sH_0(Y)\ra \sH_0(X)$ restricts to the
 desired functor
 $f_*: \sB_0^\obs(Y)\ra \sB_0^\obs(X)$.

For $f^*$,
 let $\ol x_* = (x_*, \lam_*) \in \sB_0^\obs(X)_0$
  (see Remark \ref{rk_H0_H1_struc} \eqref{it_H0_H1_struc}).
By definition, and Remark \ref{rk_A_struc} \eqref{it_0_1_ess},
 there exists some $x_0\in X(A)^\obs$ and
 $(x_*, \lam_*, \mu_*)\ra ({x_0}_*, \ep_x, \ep_x^{-1})$
 in $\sH_1(X)_1$ (who maps to  $\xi: x_*\ra {x_0}_*$ under $UV$)
 with
 $V((x_*, \lam_*, \mu_*)) = \ol x_*$, where
 $\ep_x: q_*\xra{\sim}{x_0}_*{p_X}_*$.
Using the assumption that
   $Y(A)^\obs\ra X(A)^\obs$ is  surjective,
 we know there exists some $y_0\in Y(A)^\obs$ such that $fy_0=x_0$.
Recall that
 the functor $f^*: \sH_0(X)\ra \sH_0(Y)$ sends
 $\ol x_*$ (resp.  $({x_0}_*, \ep_x)$) to
 $(f^*x_*, (\ep^{-1}f^*x_*) ({p_X}_*ax_*) \lam_*)$
 (resp.
 $(f^*{x_0}_*, (\ep^{-1}f^*{x_0}_*) ({p_X}_*a{x_0}_*) \ep_x)$), where
 $\ep: {p_Y}_*\xra{\sim}{p_X}_*f_*$,
 and $a: \id\ra f_*f^*$ is the unit map.
Let $\ep_f: {x_0}_*\xra{\sim} f_*{y_0}_*$
 be the natural equivalence
 induced by $fy_0=x_0$.
Then by relations between $a$ and $b$, and the functoriality of $\sD$,
 we see that the diagram
\ga{ \label{eq_f^*}
(f^*x_*, (\ep^{-1}f^*x_*) ({p_X}_*ax_*) \lam_*) \xra{f^*\xi_*}
 (f^*{x_0}_*, (\ep^{-1}f^*{x_0}_*) ({p_X}_*a{x_0}_*) \ep_x) \xra{f^*\ep_f}
    \\ \nonumber
 (f^*f_*{y_0}_*, ({p_Y}_*f^*\ep_f)(\ep^{-1}f^*{x_0}_*) ({p_X}_*a{x_0}_*) \ep_x)
    \xra{b{y_0}_*}
 ({y_0}_*, \ep_y)
}
 is in $\sH_0(Y)$,
  where $\ep_y: q_*\xra{\sim}{y_0}_*p_Y*$ and
  $b: f^*f_*\ra \id$ is the counit map.
Since we can lift \eqref{eq_f^*} to $\sH_1(Y)$
 (Remark \ref{rk_A_struc} \eqref{it_B0_lift}),
 it follows that $f^*\ol x_*\in \sB_0(Y)_0$ and hence
 $f^*: \sH_0(X)\ra \sH_0(Y)$ restricts to the
 desired functor
 $f^*: \sB_0^\obs(X)\ra \sB_0^\obs(Y)$.
The proof is complete.
}

\section{Cohomological descent for obstructions} \label{cohdes}

\lemm{ \label{lemm_H0_rad}
Let $f: X_0\ra X_{-1}$ be a smooth map in
 $\Chp_{/S}$ such that $X_0\in \Esp_{/S}$,  and
 $X_\bu: \N(\bDel_+)^\op \ra \N(\Chp_{/S})$ its {\Cech} nerve.
Then for any $\alpha: [m]\ra [n]$ in $\bDel_+$,
 the diagram
\eqn{
\xymatrix{
\sH_0(X_m)\ar[r]^-{d_0^*} \ar[d]^-{\alpha^*}
 &\sH_0(X_{m+1})\ar[d]^-{\alpha^*} \\
\sH_0(X_n)\ar[r]^-{d_0^*}  &\sH_0(X_{n+1}).
}}
is right adjointable.
}
\pf{
In view of the comments after \eqref{eq_D},
 in  Lemma
 \ref{lemm_D_rad},  taking $X={X_{-1}}_\liset$, and
 $U_\bu$ to be the sheaf represented by
 $X_\bu$.
Note that for $n\ge0$, $X_n\in \Esp_{/S}$.
Hence all $\alpha_*: X_n\ra X_m$ and ${d_0}_*: X_{m+1}\ra X_m$ are
 representable,
 and we may replace  $\sD(X_{/U_n})$ by $\sD(X_n)$
 \cite[Lem. 5.3.2]{lz17enhanced}, and
  obtain that
\eqn{
\xymatrix{
\sD(X_m)\ar[r]^-{d_0^*} \ar[d]^-{\alpha^*}
 &\sD(X_{m+1})\ar[d]^-{\alpha^*} \\
\sD(X_n)\ar[r]^-{d_0^*}  &\sD(X_{n+1})
}}
 is right adjointable,
 with right adjoint
\eqn{
\xymatrix{
\sD(X_m)\ar@{<-}[r]^-{{d_0}_*} \ar[d]^-{\alpha^*}
 &\sD(X_{m+1})\ar[d]^-{\alpha^*} \\
\sD(X_n)\ar@{<-}[r]^-{{d_0}_*}  &\sD(X_{n+1}).
}}
Applying $\Fun(\sD(A), -)$  and
 forming $\sH_0(-)$,  by  Lemma \ref{lemm_H0_adj},
 we obtain that
\eqn{
\xymatrix{
\sH_0(X_m)\ar[r]^-{d_0^*} \ar[d]^-{\alpha^*}
 &\sH_0(X_{m+1})\ar[d]^-{\alpha^*} \\
\sH_0(X_n)\ar[r]^-{d_0^*}  &\sH_0(X_{n+1})
}}
 is right adjointable,
 with right adjoint
\eqn{
\xymatrix{
\sH_0(X_m)\ar@{<-}[r]^-{{d_0}_*} \ar[d]^-{\alpha^*}
 &\sH_0(X_{m+1})\ar[d]^-{\alpha^*} \\
\sH_0(X_n)\ar@{<-}[r]^-{{d_0}_*}  &\sH_0(X_{n+1}).
}}
}

\lemm{ \label{lemm_H0_lim}
Let  $\La\in\Ring_{\Box\tu{-tor}}$.
Let $f: X_0\ra X_{-1}$ be a smooth
 surjective map in
 $\Chp_{/S}^\Box$ such that $X_0\in \Esp_{/S}$,  and
 $X_\bu: \N(\bDel_{s,+})^\op \ra \N(\Chp_{/S}^\Box)$
 its semisimplicial {\Cech} nerve $($see  \cite[Notation 6.5.3.6]{HTT}$)$.
Then there is an canonical equivalence
\eq{ \label{eq_H0_lim}
\sH_0(X_{-1})\ra \varprojlim_{n\in\bDel_s}\sH_0(X_n),
}
 where the maps are using  $*$-pullback.
}
\pf{
Extend $X_\bu$  to a full {\Cech} nerve
 $X_\bu: \N(\bDel_+)^\op \ra \N(\Chp_{/S}^\Box)$.
We first show that
 there is an canonical equivalence
\eq{ \label{eq_H0_lim_full}
\sH_0(X_{-1})\ra \varprojlim_{n\in\bDel}\sH_0(X_n),
}
 where the maps are using  $*$-pullback.
We apply \cite[Lem. 3.3.6]{lz17enhanced}.
For  every $X\in \Chp_{/S}$,
By  Lemma \ref{lemm_H0_U_lim} \eqref{it_H0_adm_lim},
 $\sH_0(X)$  admits
 limits of  cosimplicial  objects.
Since
  $\La\in\Ring_{\Box\tu{-tor}}$ (recall that we write  $\sD(-)=\sD(-,\La)$)
  and
  $f: X_0\ra X_{-1}$ is a smooth
 surjective map in
 $\Chp_{/S}^\Box$,
 by {\Poin} duality (\cite[Thm. 6.2.9]{lz17enhanced}),
 $f^*\xra{\sim} \langle-\dim f\rangle f^!$
 has a left adjoint $\langle\dim f\rangle f_!$.
Thus  $f^*: \sH(X_{-1})\ra \sH(X_0)$  preserves small limits.
In view of  the homotopy commutative diagram
\eq{\label{eq_U_f^*}
\xymatrix{
\sH_0(X_{-1})\ar[r]^-{U}\ar[d]^-{f^*} &\sH(X_{-1})\ar[d]^-{f^*} \\
\sH_0(X_0)\ar[r]^-{U} &\sH(X_0),
}}
 by  Lemma \ref{lemm_H0_U_lim} \eqref{it_U},
 the restriction $f^*: \sH_0(X_{-1})\ra \sH_0(X_0)$ also
 preserves small limits.
Thus assumption (1) follows.
Assumption (2) follows from Lemma \ref{lemm_H0_rad}.
For Assumption (3),  first note that $U$ is conservative
  by  Lemma \ref{lemm_H0_U_lim} \eqref{it_U} (since
 $x\ra y$ is an equivalence if and only if it is a limit diagram).
Then  by \cite[Lem. 4.3.4]{lz17enhanced}, $f^*$ is conservative
 since $f$ is smooth surjective.
Use \eqref{eq_U_f^*} again,
 we obtain that $f^* \sH_0(X_{-1})\ra \sH_0(X)$ is conservative,
 i.e., assumption (3) follows.
It follows that \eqref{eq_H0_lim_full} holds.
Then by the dual version of \cite[Lem. 6.5.3.7]{HTT},
 the diagram in  \eqref{eq_H0_lim_full} can be restricted
 to $\bDel_s$, and resulting in the equivalence \eqref{eq_H0_lim}.
The proof is complete.
}
%

\cor{ \label{cor_B0_lim_weak}
Let  $\La\in\Ring_{\Box\tu{-tor}}$.
Let $f: X_0\ra X_{-1}$ be a smooth surjective
 map in $(\Chp_{/S}^\sFoB)_{\sFo}$
 such that $X_0\in \Esp_{/S}^\sFoB$,  and
\eq{ \label{eq_sCech}
X_\bu: \N(\bDel_{s,+})^\op \ra \N(\Chp_{/S}^\sFoB)_{\sFo}
}
 its semisimplicial {\Cech} nerve.
Then there is an canonical equivalence
\eqn{
\sB_0^\obs(X_{-1})\ra \varprojlim_{n\in\bDel_s}\sB_0^\obs(X_n),
}
 where the maps are using  $*$-pullback.
}
\rk{
First we see that  \eqref{eq_sCech} makes sense
 since $\sFo$ is stable under composition and pullback.
Then   by Lemma \ref{lemm_B0_adj},
 the diagram $\sB_0^\obs\circ X_\bu^\op:  N(\bDel_s) \ra \Cat_\infty$
 using $*$-pullback
 also makes sense since we can do $*$-pullback along maps in $\sFo$.
}

\pf{
Note that since $\sB_0^\obs(X)$ is the
  essential image of $V$, it is strictly full subcategory of $\sH_0(X)$
  for any $X$.
We apply  Lemma \ref{lemm_H0_lim} and
 \cite[Lem. 3.1.4]{lz17enhanced}.
It suffices to show that   for any $\ol x_*\in \sH_0(X_{-1})$,
 $f^*\ol x_*\in \sB_0^\obs(X_0)$ implies $\ol x_*\in \sB_0^\obs(X_{-1})$.
Indeed,
 in view of Lemmas \ref{lemm_H0_adj} and \ref{lemm_B0_adj},
 apply $f_*$ to  $f^*\ol x_*\in \sB_0^\obs(X_0)$,
  and we obtain   the unit map $\ol x_* \ra f_*f^*\ol x_*$
 in $\sH_0(X_{-1})$
 with
 $f_*f^*\ol x_*\in \sB_0^\obs(X_{-1})$.
By Remark \ref{rk_A_struc} \eqref{it_B0_lift},
 we obtain that
 $\ol x_* \in  \sB_0^\obs(X_{-1})$.
The proof is complete.
}

%

Next we further
 remove the assumption that
 $X_0\in \Esp_{/S}^\sFoB$
  in Corollary \ref{cor_B0_lim_weak}.

\thm{ \label{thm_B0_lim}
Let  $\La\in\Ring_{\Box\tu{-tor}}$.
Let $\sE$ be a  family of maps of $\Chp_{/S}$  such that
\enmt{[\upshape (a)]
\item $\sE$ contains every degenerate edge,
\item $\sE$  are stable under composition and pullback,
\item every $f\in \sE$ is smooth surjective, and
\item $\sE\subseteq \sFoB$.
}
Then $\sB_0^\obs$ defines a functor
\eq{ \label{eq_B0_fun}
\sB_0^\obs: \N(\Chp_{/S}^{\sE})_{\sFo}^\op \ra \Cat_\infty
}
 such that for any
 $f: X_0\ra X_{-1}$ being a map
 in $(\Chp_{/S}^{\sE})_{\sE}$, and
 $X_\bu: \N(\bDel_{s,+})^\op \ra \N(\Chp_{/S}^{\sE})_{\sE}$
 its semisimplicial {\Cech} nerve,
 the map is an equivalence
\eq{ \label{eq_B0_lim}
\sB_0^\obs(X_{-1})\ra \varprojlim_{n\in\bDel_s}\sB_0^\obs(X_n).
}
}
\pf{
We first show the result for  the category $\Esp_{/S}^\sE$.
Let $f: X_0\ra X_{-1}$ be a map
 in $(\Esp_{/S}^\sE)_{\sE}$, and
 $X_\bu: \N(\bDel_{s,+})^\op \ra \N(\Esp_{/S}^\sE)_{\sE}$
 its semisimplicial {\Cech} nerve.
With assumptions  of Corollary \ref{cor_B0_lim_weak}
 fulfilled, and using Lemma \ref{lemm_B0_adj},
 the result for  the category $\Esp_{/S}^\Box$ holds, that is,
 $\sB_0^\obs$ defines a functor
\eq{ \label{eq_B0_fun_Esp}
\sB_0^\obs: \N(\Esp_{/S}^\sE)_{\sFo}^\op \ra \Cat_\infty
}
 via $*$-pullback,
 such that
\eq{ \label{eq_B0_lim_Esp}
\sB_0^\obs(X_{-1})\xra{\sim} \varprojlim_{n\in\bDel_s}\sB_0^\obs(X_n).
}
 $f: X_0\ra X_{-1}$ being a map
 in $(\Esp_{/S}^\sE)_{\sE}$.

Now we consider  the category  $\Chp_{/S}^\sE$.
Fix $X_{-1}\in  \Chp_{/S}^\sE$ and choose an $\sE$-atlas
 $f: X_0\ra X_{-1}$, that is, $f$ is  a  map in
 $\sE$ such that $X_0\in \Esp_{/S}^\sE$.
Then Proposition \ref{cor_B0_lim_weak}  is still available,
 showing  that
\eq{ \label{eq_B0_lim_atlas}
\text{\eqref{eq_B0_lim} still holds for the above $f$.}
}

Now we consider the functoriality of $\sB_0^\obs$, i.e.,
 to show it is a functor \eqref{eq_B0_fun}.
We apply \cite[Lem. 5.3.1]{lz17enhanced} with $\sC = \N(\Esp_{/S}^\sE)$,
 $\tilde \sC = \N(\Chp_{/S}^\sE)$, $\tilde\sE = \sE$,
 $\tilde\sF = \sFo$
 and $\sD = \Cat_\infty$ (by
  \cite[Cor. 4.2.4.8]{HTT} $\Cat_\infty$ admits totalizations).
It follows that the restriction
\eq{ \label{eq_Fun_des_ext}
\Fun^\sE(\N(\Chp_{/S}^\sE)_\sFo^\op, \Cat_\infty)\ra
 \Fun^\sE(\N(\Esp_{/S}^\sE)_\sFo^\op, \Cat_\infty)
}
 is a trivial fibration.
Then since the functor \eqref{eq_B0_fun_Esp}  satisfies
 \eqref{eq_B0_lim_Esp}
 for
 all $f$ in  $(\Esp_{/S}^\sE)_\sE$,
 it is in
 $\Fun^\sE(\N(\Esp_{/S}^\sE)_\sFo^\op, \Cat_\infty)$.
By \eqref{eq_Fun_des_ext}, it  extends uniquely to a functor
\eqn{
\tilde\sB_0^\obs: \N(\Chp_{/S}^\sE)_\sFo^\op\ra \Cat_\infty
}
 in
 $\Fun^\sE(\N(\Chp_{/S}^\sE)_\sFo^\op, \Cat_\infty)$.
Thus  $\tilde\sB_0^\obs$ satisfies \eqref{eq_B0_lim} for
 all  $f$ in  $(\Chp_{/S}^\sE)_\sE$.
Along with \eqref{eq_B0_lim_atlas},
 we obtain that $\tilde\sB_0^\obs$ coincides with $\sB_0^\obs$ on vertices.
It follows that  $\sB_0^\obs$ defines a functor
\eqn{
\sB_0^\obs: \N(\Chp_{/S}^\sE)_\sFo^\op\ra \Cat_\infty
}
 such that  \eqref{eq_B0_lim} is an equivalence
 for $f: X_0\ra X_{-1}$ being a  map in
 $(\Chp_{/S}^\sE)_\sE$ (since $\tilde\sB_0^\obs$ is such a functor).
The proof is complete.
}

\rk{ \label{rk_thm_B0_lim}
By Lemma \ref{lemm_B0_adj} we know that
 $\sB_0^\obs$ defines a functor
\eqn{
\sB_0^\obs: \N(\Chp_{/S}^\sE)_\sFo^\op\ra \Cat_\infty
}
 via $*$-pullback.
By Corollary \ref{cor_B0_lim_weak} we know that
 for $f: X_0\ra X_{-1}$ being a
 map in $(\Chp_{/S}^\sE)_{\sE}$
 such that $X_0\in \Esp_{/S}^\sE$,  and
\eqn{
X_\bu: \N(\bDel_{s,+})^\op \ra \N(\Chp_{/S}^\sE)_{\sE}
}
 its semisimplicial {\Cech} nerve,
 there is an canonical equivalence
\eqn{
\sB_0^\obs(X_{-1})\ra \varprojlim_{n\in\bDel_s}\sB_0^\obs(X_n).
}
Then  for general
  $f: X_0\ra X_{-1}$ being
 a map in
 $(\Chp_{/S}^\sE)_\sE$,
  we may apply
  \cite[Lem. 3.1.2 (3) and (4)]{lz17enhanced} on a chart of $f$
 to obtain  \eqref{eq_B0_lim}.
This gives another proof of Theorem \ref{thm_B0_lim}.
}

\rk{ \label{rk_B0_sub}
\enmt{[\upshape (1)]
\item \label{it_B0_edge}
For an map $f: Y\ra X$ in  $(\Chp_{/S}^\sE)_\sFo$,
 we denote by $f^\natural: \sB_0^\obs(X)\ra \sB_0^\obs(Y)$
 the image of $f$ under \eqref{eq_B0_fun}.
By Corollary \ref{cor_B0_lim_weak} and Lemma \ref{lemm_B0_adj},
 if $f$ is in $(\Esp_{/S}^\sE)_\sFo$,
 then $f^\natural$ is equivalent to $f^*$.
Then by the functoriality  of  the construction of
 $\tilde \sB_0^\obs$, we deduce that for  arbitrary
 $f: Y\ra X$ in  $(\Chp_{/S}^\sE)_\sFo$,
  $f^\natural$ is equivalent to $f^*$.

Moreover,
 by Lemma \ref{lemm_B0_adj},
 we have
 $\sB_0^\obs\in \Fun^\RAd(\N(\Chp_{/S}^{\sE})_{\sFo}^\op,  \Cat_\infty)_0$
 with right adjoint given by $*$-pushforward.

\item \label{it_B0_H}
Use Lemma \ref{lemm_H0_adj} instead of Lemma \ref{lemm_B0_adj}.
Then we see that in
 Theorem \ref{thm_B0_lim},
 we may replace $\sFo$ by $\rep$ and $\sB_0^\obs$ by $\sH_0$.
To include this case, we use notation $\sFm=\rep$ and $\sB_0^m = \sH_0$
 in
 Theorem \ref{thm_B0_lim},
}}

\constr{ \label{constr_C}
Let  $\La\in\Ring_{\Box\tu{-tor}}$.
Suppose that for every $X\in \Chp_{/S}$, we have a subset
 $X(A)^\olobs \subseteq X(A)$,
 and for every $X\in \Sch_{/S}$, we have a subset
\eqn{
X(A)^\obs\subseteq   X(A)^\olobs \subseteq X(A).
}
Suppose that  $\sE$ and $\sFolo$ satisfy the conditions for $\sE$ and $\sFo$
 in Theorem \ref{thm_B0_lim}.
By  Theorem \ref{thm_B0_lim} and Remark \ref{rk_B0_sub} \eqref{it_B0_edge},
 we know that $\sB_0^\olobs\in
 \Fun^{\RAd, \sE}(\N(\Sch_{/S})_\sFolo^\op, \Cat_\infty)_0$.
Consider the full subcategory
\eqn{
\sK \subseteq
 \Fun^{\RAd, \sE}(\N(\Sch_{/S})_\sFolo^\op, \Cat_\infty)_{/\sB_0^\olobs}
}
 spanned by  $\sC_0\ra \sB_0^\olobs$ such that
 there is two fully faithful embeddings of functors
 $\sA_0^\obs(X)\hra \sC_0(X) \hra \sB_0^\olobs(X)$
 for every $X\in \Sch_{/S}$, and that  $\sC_0(X)$ is strict when
 viewed as a fully faithful subcategory of $\sB_0^\olobs(X)$.
Clearly $\sK\neq \emptyset$ since it has $\sB_0^\olobs\xra{\id} \sB_0^\olobs$
 as a
 final object.
Since
 $\Fun^\RAd(\N(\Sch_{/S})_\sFolo^\op, \Cat_\infty)$
 admits small limits  \cite[Cor. 4.7.4.18 (1)]{HA},
 and one may swap two compatible limits,
 along with \cite[Cor. 5.1.2.3]{HTT},
 we know that
 the category
 $\Fun^{\RAd, \sE}(\N(\Sch_{/S})_\sFolo^\op, \Cat_\infty)_{/\sB_0^\olobs}$
 also admits small limits.
Then we denote by $\sC^\obsol\ra \sB_0^\olobs$  the limit
\eq{ \label{eq_C}
\varprojlim_{\sK}(\sC_0\ra \sB_0^\olobs)  \in
 \Fun^{\RAd, \sE}(\N(\Sch_{/S})_\sFolo^\op, \Cat_\infty)_{/\sB_0^\olobs},
}
 whose source $\sC^\obsol\in
 \Fun^{\RAd, \sE}(\N(\Sch_{/S})_\sFolo^\op, \Cat_\infty)_0$.

Extend $\sC^\obsol$ by the trivial fibrations
\eq{ \label{eq_Fun_des_exts}
\Fun^\sE(\N(\Chp_{/S}^\sE)_\sFolo^\op, \Cat_\infty)\ra
 \Fun^\sE(\N(\Esp_{/S}^\sE)_\sFolo^\op, \Cat_\infty)\ra
 \Fun^\sE(\N(\Sch_{/S})_\sFolo^\op, \Cat_\infty)
}
 to obtain a functor (unique up to equivalence)
\eq{ \label{eq_C_fun}
\sC^\obsol: \N(\Chp_{/S}^\sE)_\sFolo^\op\ra \Cat_\infty
}
 such that for any
 $f: X_0\ra X_{-1}$ being a map
 in $(\Chp_{/S}^{\sE})_{\sE}$, and
 $X_\bu: \N(\bDel_{s,+})^\op \ra \N(\Chp_{/S}^{\sE})_{\sE}$
 its semisimplicial {\Cech} nerve,
 the map is an equivalence
\eq{ \label{eq_C_lim}
\sC^\obsol(X_{-1})\ra \varprojlim_{n\in\bDel_s}\sC^\obsol(X_n).
}
}

\defi{ \label{defi_inc_fun}
Let $W$ be a simplicial set.
We say an edge $F\ra G$ in  $\Fun(W, \Cat_\infty)$
 is an \emph{inclusion of functors} if
 for each $w\in W$, the functor
 $F(w)\ra G(w)$ is fully faithful,
 and   $F(w)$ is strict when viewed as a full subcategory of $G(w)$.
For this we use the notation $F\subseteq G$.
}

\rk{ \label{rk_C}
Let
 $f: X_0\ra X_{-1}$ be a map
 in $(\Chp_{/S}^{\sE})_{\sE}$, and
 $X_\bu: \N(\bDel_{s,+})^\op \ra \N(\Chp_{/S}^{\sE})_{\sE}$
 its semisimplicial {\Cech} nerve.
\enmt{[\upshape (1)]
\item \label{it_C_sub}
By the above construction \eqref{eq_C}, \eqref{eq_Fun_des_exts},
 and Theorem \ref{thm_B0_lim}, Remark \ref{rk_B0_sub} \eqref{it_B0_edge},
 \eqref{eq_C_lim}
 we have the inclusion of functors
\eqn{
\sC^\obsol \subseteq \sB_0^\olobs: \N(\Chp_{/S}^\sE)_\sFolo^\op\ra \Cat_\infty
}
 such that
\eqn{ 
\xymatrix{
\sC^\obsol(X_{-1})\ar[r]^-{\sim}\ar@{^(->}[d]
 &\varprojlim_{n\in\bDel_s}\sC^\obsol(X_n)
 \ar@{^(->}[d]  \\
\sB_0^\olobs(X_{-1})\ar[r]^-{\sim}
 &\varprojlim_{n\in\bDel_s}\sB_0^\olobs(X_n)
}}
 commutes up to homotopy,
 and one has a fully faithful embedding of functors
 $\sA^\obs(X)\hra\sC^\obsol(X)$
 for every $X\in \Sch_{/S}$.
In particular, the functor \eqref{eq_C_fun} on edges
 agrees with $*$-pullback.

Moreover, by and
  \eqref{eq_C}, the limit construction of the
   extensions \eqref{eq_Fun_des_exts},
  and
 \cite[Lem. 4.3.7]{lz17enhanced},
 we have that $\sC^\obsol\ra \sB_0^\olobs$ is an edge of
 $\Fun^{\RAd, \sE}(\N(\Chp_{/S}^{\sE})_{\sFolo}^\op,  \Cat_\infty)$.
In particular, by  Remark \ref{rk_B0_sub} \eqref{it_B0_edge},
 $\sC^\obsol\in
  \Fun^{\RAd, \sE}(\N(\Chp_{/S}^{\sE})_{\sFolo}^\op,  \Cat_\infty)_0$
 with right adjoint given by $*$-pushforward.

\item  \label{it_C_cmp}
Suppose that for all $X\in \Sch_{/S}$, we are given  subsets
 $X(A)^{\obs_1}\subseteq X(A)^{\obs_2} \subseteq X(A)^\olobs\subseteq X(A)$.
Then by definition
 $\sA^{\obs_1}(X)\subseteq \sA^{\obs_2}(X)$ for all $X\in \Sch_{S}$.
It follows that
 we have the inclusion of functors
\eqn{
\sC^{\obs_1,\olobs} \subseteq \sC^{\obs_2,\olobs}:
 \N(\Sch_{/S})_\sFolo^\op\ra \Cat_\infty
}
Then by \eqref{eq_Fun_des_exts} again, the inclusion extends to
\eqn{
\sC^{\obs_1,\olobs} \subseteq \sC^{\obs_2,\olobs}:
 \N(\Chp_{/S}^\sE)_\sFolo^\op\ra \Cat_\infty
}
 such that
\eqn{
\xymatrix{
\sC^{\obs_1,\olobs}(X_{-1})\ar[r]^-{\sim}\ar@{^(->}[d]
 &\varprojlim_{n\in\bDel_s}\sC^{\obs_1,\olobs}(X_n)
 \ar@{^(->}[d]  \\
\sC^{\obs_2,\olobs}(X_{-1})\ar[r]^-{\sim}
 &\varprojlim_{n\in\bDel_s}\sC^{\obs_2,\olobs}(X_n)
}}
 commutes up to homotopy.

\item \label{it_C_E}
Let $\sE_0\subset\sE$ be another family satisfying conditions
 in Theorem \ref{thm_B0_lim}.
Then  \eqref{eq_Fun_des_exts} fits into
 the following homotopy commutative diagram
\eq{ \label{eq_res}
\xymatrix{
\Fun^\sE(\N(\Chp_{/S}^\sE)_\sFolo^\op, \Cat_\infty)\ar[r]^-{\res_c}
  \ar[d]^-{\res_1}
 &\Fun^{\sE_0}(\N(\Chp_{/S}^{\sE_0})_\sFolo^\op, \Cat_\infty)
  \ar[d]^-{\res_{1,0}} \\
\Fun^\sE(\N(\Esp_{/S}^\sE)_\sFolo^\op, \Cat_\infty)\ar[r]^-{\res_e}
  \ar[d]^-{\res_2}
 & \Fun^{\sE_0}(\N(\Esp_{/S}^{\sE_0})_\sFolo^\op, \Cat_\infty)
  \ar[d]^-{\res_{2,0}} \\
\Fun^\sE(\N(\Sch_{/S})_\sFolo^\op, \Cat_\infty)\ar@{^(->}[r]^-{i_s}
 &\Fun^{\sE_0}(\N(\Sch_{/S})_\sFolo^\op, \Cat_\infty) \\
\Fun^{\RAd, \sE}(\N(\Sch_{/S})_\sFolo^\op, \Cat_\infty)_{/\sB_0^\olobs}
  \ar[r]^-{i_b}\ar[u]
 &\Fun^{\RAd, {\sE_0}}(\N(\Sch_{/S})_\sFolo^\op, \Cat_\infty)_{/\sB_0^\olobs}
  \ar[u] \\
\sK\ar[r]^-{i_k}\ar@{^(->}[u] &\sK_0\ar@{^(->}[u]
}}
 where the restrictions $\res_1$, $\res_2$, $\res_{1,0}$ and $\res_{2,0}$
 are trivial fibrations, and $i_s$ is a fully faithful embedding which induces
 $i_b$ and $i_k$.
Thus the maps $\res_c$ and $\res_e$ are also fully faithful embeddings,
 and are
 restrictions along the fully faithful embeddings
 $\N(\Chp_{/S}^{\sE_0})_\sFolo^\op\ra \N(\Chp_{/S}^\sE)_\sFolo^\op$ and
 $\N(\Esp_{/S}^{\sE_0})_\sFolo^\op\ra \N(\Esp_{/S}^\sE)_\sFolo^\op$,
 respectively.
One checks that all horizontal arrows preserve small limits.

Denote  by
\eqn{
\sC_0^\obsol: \N(\Chp_{/S}^{\sE_0})_\sFolo^\op\ra \Cat_\infty
}
 the functor in Construction \ref{constr_C} corresponding to $\sE_0$.
It follows that  we have   inclusions of functors
\eqn{
\sC_0^\obsol \subseteq \res_c\sC^\obsol \subseteq \sB_0^\olobs:
 \N(\Chp_{/S}^{\sE_0})_\sFolo^\op\ra \Cat_\infty.
}
}}

\defi{ \label{defi_obs_functorial}
Let $\Chp_{/S}'\subseteq \Chp_{/S}$ be a sub $2$-category.
Suppose that attached to each $X\in \Chp_{/S}'$,
 we have a subset $X(A)^\obs\subseteq X(A)$.
We say that $\obs$ is \emph{functorial} on $\Chp_{/S}'$
 if the assignment $X\mpt X(A)^\obs$ gives a functor
\eqn{
-(X)^\obs: \Chp_{/S}'\ra \Set.
}}

\rk{ \label{rk_functorial}
In particular, if
 $\obs: (\Chp_{/S}')^\op\ra \Set$
 is a stable functor, and
 for $X\in \Chp_{/S}'$, let
 $X(A)^\obs$ the corresponding
 obstruction set (see Section \ref{ptob}),
 then  $\obs$ is functorial on  $\Chp_{/S}'$.
}

\thm{ \label{thm_C}
Let  $\La\in\Ring_{\Box\tu{-tor}}$.
Suppose that for every $X\in \Chp_{/S}$, we have a subset
 $X(A)^\olobs \subseteq X(A)$,
 and for every $X\in \Sch_{/S}$, we have a subset
\eqn{
X(A)^\obs\subseteq   X(A)^\olobs \subseteq X(A).
}
Fix an $\sFolo$.
Let $\sE$ be a  family of maps of $\Chp_{/S}$  such that
\enmt{[\upshape (a)]
\item $\sE$ contains every degenerate edge,
\item $\sE$  are stable under composition and pullback,
\item every $f\in \sE$ is smooth surjective, and
\item $\sE\subseteq \sFolo\cap \Box$.
}
Let
\eqn{
\obs: \Sch_{/S}\ra \Set, \quad X\mpt X(A)^\obs\subseteq X(A)^\olobs
}
 be a map.
Then There is a way to associate to
 every pair $(\obs, \sE)$
 a functor
\eqn{
\sC_\sE^\obsol: \N(\Chp_{/S}^\sE)_\sFolo^\op\ra \Cat_\infty
}
 satisfying the following properties.
\enmt{[\upshape (1)]
\item \label{it_C_A}
For any $X\in \Sch_{/S}$,
 we have a fully faithful  functor
  $\sA_0^\obs(X)\hra \sC_\sE^\obsol(X)$.

\item \label{it_C_des}
For any
 $f: X_0\ra X_{-1}$ being a map
 in $(\Chp_{/S}^{\sE})_{\sE}$, and
 $X_\bu: \N(\bDel_{s,+})^\op \ra \N(\Chp_{/S}^{\sE})_{\sE}$
 its semisimplicial {\Cech} nerve,
 the map is an equivalence
\eqn{
\sC_\sE^\obsol(X_{-1})\ra \varprojlim_{n\in\bDel_s}\sC_\sE^\obsol(X_n).
}

\item \label{it_C_sub1}
We have  that
 $\sC_\sE^\obsol \ra \sB_0^\olobs$ is an edge of
 $\Fun^{\RAd, \sE}(\N(\Chp_{/S}^{\sE})_{\sFolo}^\op,  \Cat_\infty)$
 making
 an inclusion of functors
\eqn{
\sC_\sE^\obsol \subseteq \sB_0^\olobs:
 \N(\Chp_{/S}^\sE)_\sFolo^\op\ra \Cat_\infty
}
 such that
\eqn{
\xymatrix{
\sC_\sE^\obsol(X_{-1})\ar[r]^-{\sim}\ar@{^(->}[d]
 &\varprojlim_{n\in\bDel_s}\sC_\sE^\obsol(X_n)
 \ar@{^(->}[d]  \\
\sB_0^\olobs(X_{-1})\ar[r]^-{\sim}
 &\varprojlim_{n\in\bDel_s}\sB_0^\olobs(X_n)
}}
 commutes up to homotopy for
 $f: X_0\ra X_{-1}$ being a map as in \eqref{it_C_des}.
In particular, the functor $\sC_\sE^\obsol$ on edges
 agrees with $*$-pullback, with right adjoint agrees with $*$-pushforward.

\item \label{it_C_cmp1}
Suppose that we are given two maps
\eqn{
\obs_1\subseteq \obs_2: \Sch_{/S}\ra \Set
}
Then we have an inclusion of functors
\eqn{
\sC_\sE^{\obs_1,\olobs} \subseteq \sC^{\obs_2,\olobs}:
 \N(\Chp_{/S}^\sE)_\sFolo^\op\ra \Cat_\infty
}
 such that
\eqn{
\xymatrix{
\sC_\sE^{\obs_1,\olobs}(X_{-1})\ar[r]^-{\sim}\ar@{^(->}[d]
 &\varprojlim_{n\in\bDel_s}\sC_\sE^{\obs_1,\olobs}(X_n)
 \ar@{^(->}[d]  \\
\sC^{\obs_2,\olobs}(X_{-1})\ar[r]^-{\sim}
 &\varprojlim_{n\in\bDel_s}\sC^{\obs_2,\olobs}(X_n)
}}
 commutes up to homotopy for
 $f: X_0\ra X_{-1}$ being a map as in \eqref{it_C_des}.

\item \label{it_C_B}
If
 $\obs$ extends  to a map
\eqn{
\obs: \Chp_{/S}^\sE\ra \Set, \quad X\mpt X(A)^\obs\subseteq X(A)
}
 and is functorial on  $\Chp_{/S}^\sE$,
 and that for each $f: Y\ra X$ in $\sE\cap (\Sch_{/S})_1$,
  $f$  induces a surjective map
  $Y(A)^\obs\ra X(A)^\obs$,
 then we have a fully faithful  functor
 $\sC_\sE^\obsol(X)\hra \sB_0^\obs(X)$ for every $X\in\Chp_{/S}^\sE$.

\item \label{it_C_E1}
Let $\sE_0\subset\sE$ be another family satisfying conditions
 at the begin of the theorem.
Then we have   inclusions of functors
\eqn{
\sC_{\sE_0}^\obsol \subseteq \res_c\sC_\sE^\obsol
 \subseteq \sB_0^\olobs:
 \N(\Chp_{/S}^{\sE_0})_\sFolo^\op\ra \Cat_\infty
}
 compatible with the equivalences in \eqref{it_C_des}, i.e.,
\eqn{
\xymatrix{
\sC_{\sE_0}^\obsol(X_{-1})\ar[r]^-{\sim}\ar@{^(->}[d]
 &\varprojlim_{n\in\bDel_s}\sC_{\sE_0}^\obsol(X_n)
 \ar@{^(->}[d]  \\
\res_c\sC_\sE^\obsol(X_{-1})\ar[r]^-{\sim}
 &\varprojlim_{n\in\bDel_s}\res_c\sC_\sE^\obsol(X_n)
}}
 commutes up to homotopy for
 $f: X_0\ra X_{-1}$ being a map as in \eqref{it_C_des},
 where $\res_c$ is induced by the fully faithful embedding
 $\N(\Chp_{/S}^{\sE_0})_\sFolo^\op\hra
  \N(\Chp_{/S}^\sE)_\sFolo^\op$ $($see \eqref{eq_res}$)$.
}
Moreover, for fixed $(\obs, \sE)$, the functor
 $\sC_\sE^\obsol$, up to equivalence, is the smallest one satisfying
 \eqref{it_C_A}, \eqref{it_C_des} and \eqref{it_C_sub1},
 that is, for any  ${\sC'}_\sE^\obsol$ satisfying
 \eqref{it_C_A}, \eqref{it_C_des} and \eqref{it_C_sub1},
 we have an inclusion of functors  ${\sC'}_\sE^\obsol\subseteq
  \sC_\sE^\obsol$.
}

\pf{
We use Construction \ref{constr_C}.
Then  the resulting $\sC^\obsol$ is a candidate for
 $\sC_\sE^\obsol$.
It suffices to verify the properties.
Indeed,  \eqref{it_C_A} and \eqref{it_C_des}
  follows from the construction
  of $\sC^\obsol$  \eqref{eq_C} and \eqref{eq_C_lim}.
The property \eqref{it_C_sub1} following from  Remark
 \ref{rk_C} \eqref{it_C_sub}.
The smallestness also follows from
 \eqref{eq_C} and \eqref{eq_Fun_des_exts}.
The property \eqref{it_C_cmp1} following from  Remark
 \ref{rk_C} \eqref{it_C_cmp}.

To show \eqref{it_C_B},
 note that by
 the assumption
 that for each $f: Y\ra X$ in $\sE\cap (\Sch_{/S})_1$,
  $f$  induces a surjective map
  $Y(A)^\obs\ra X(A)^\obs$,
 we may apply
   Theorem \ref{thm_B0_lim} on $\Sch_{/S}$ with $\sE=\sFo$
   being $\sE\cap (\Sch_{/S})_1$,
 to obtain that
 $\sB_0^\obs$ is in
 $\Fun^\sE(\N(\Sch_{/S})_\sFolo^\op, \Cat_\infty)$.
It follows by \eqref{eq_C} that we have an inclusion of functors
\eq{ \label{eq_C_B_Sch}
\sC_\sE^\obsol\subseteq \sB_0^\obs:
  \N(\Sch_{/S})_\sFolo^\op\ra \Cat_\infty.
}
Let  $X\in \Chp_{/S}^\sE$ and
 $x\in \sC_\sE^\obsol(X)$.
Choose an $\sE$-atlas  for $X$,
 that is  a smooth surjective map
 $f: X_0\ra X_{-1} = X$
 in $(\Chp_{/S}^{\sE})_{\sE}$ with $X_0\in \Esp_{S}^\sE$.
By \eqref{it_C_sub1}, the functor
\eqn{
\sC_\sE^\obsol: \N(\Chp_{/S}^\sE)_\sFolo^\op\ra \Cat_\infty,
}
 on edges
 agrees with $*$-pullback, with right adjoint agrees with $*$-pushforward.
Thus we know that
 $f^*x\in  \sC_\sE^\obsol(X_0)$.
Since
 $\obs$ is functorial on $\Chp_{/S}^\sE$,
  one checks that
 $f_*: \sH_0(X_0)\ra \sH_0(X_{-1})$ sends
  $\sB_0^\obs(X_0)$ to $\sB_0^\obs(X_{-1})$.
Thus if  \eqref{it_C_B} holds for  $X_0$,
 then we have the unit map $x_* \ra f_*f^*x_*$
 in $\sH_0(X_{-1})$
 with
 $f_*f^*x_*\in \sB_0^\obs(X_{-1})$.
Again use
 Remark \ref{rk_A_struc} \eqref{it_B0_lift}
 to obtain that
 $x_* \in  \sB_0^\obs(X_{-1})$.
It follows that we reduce to show \eqref{it_C_B} for
  $X\in \Esp_{/S}^\sE$, and then  for
  $X\in \Sch_{/S}$.
This is true by \eqref{eq_C_B_Sch} and then
 \eqref{it_C_B} follows.

Noting that $\res_c$ preserves small limits,
 the property  \eqref{it_C_E1} directly follows from
 Remark
 \ref{rk_C} \eqref{it_C_E}.
The proof is complete.
}

\rk{ \label{rk_Sch_sub}
We may  replace $\Sch_{/S}$ by a full subcategory
 $\Sch_{/S}'$, and
 $\Esp_{/S}^\sE$ (resp. $\Chp_{/S}^\sE$)
 by the corresponding full subcategory (resp. sub $2$-category)
 ${\Esp_{/S}'}^\sE$ (reps. ${\Chp_{/S}'}^\sE$),
 that is, the full subcategory (resp. sub $2$-category) of
 ${\Esp_{/S}'}^\sE$ (reps. ${\Chp_{/S}'}^\sE$) spanned by objects that
 admits $\sE$-atlas from
 $\Sch_{/S}'$ (resp.  ${\Esp_{/S}'}^\sE$).
Then Theorem \ref{thm_C} is still correct.
To see this, just start over from
 $\Sch_{/S}'$ in Construction \ref{constr_C}, and the argument
 is the same.
}

\section{Application to obstructions on algebraic stacks} \label{appl}

We keep
 the notations and conditions in Theorem \ref{thm_C}.
From now on, we take $\olobs=m$, such that
 $\sF^\olobs=\rep$ and $\sB_0^\olobs = \sH_0$ (see Remark
 \ref{rk_B0_sub} \eqref{it_B0_H}) and use notations and
 result from Section \ref{cohdes}

Let  $X\in \Chp_{/S}^\sE$.
Recall that for $x\in X(A)$,
 $\ol x_* = (x_*, \ep)\in \sH_0(X)_0$
  where $\ep: q_*\xra{\sim}x_*p_*$ is
  the natural equivalence.

Take $\olobs = m$ so that $\sFm=\rep$ and $\sB_0^m = \sH_0$
 (see Remark \ref{rk_B0_sub} \eqref{it_B0_H}.
We recall that  Theorem \ref{thm_C} produces a functor
\eqn{
\sC_\sE^\obsm: \N(\Chp_{/S}^\sE)_\rep^\op\ra \Cat_\infty.
}

\defi{ \label{defi_tobs}
Let  $X\in \Chp_{/S}^\sE$.
We define the subset
\gan{
X(A)^{\tobs_\sE} = X(S)\ \bigcup\
 \{x\in X(A) \mid \text{ the object
  $\ol x_* = (x_*, \ep)\in \sH_0(X)_0$ } \\
   \text{ lies in the essential image of $\sC_\sE^\obsm(X)\hra \sH_0(X)$} \}
   \subseteq X(A).
}}
\rk{ \label{rk_tobs}
\enmt{[\upshape (1)]
\item By Theorem \ref{thm_C} \eqref{it_C_sub1} and the definition
 of $\sB_0$,
 we have  inclusions of functors
\eqn{
\sC_\sE^\obsm \subseteq \sH_0:
 \N(\Chp_{/S}^\sE)_\rep^\op\ra \Cat_\infty,
}
 the Definition \ref{defi_tobs} makes sense,
 and the essential image can be replaced by image.
\item \label{it_Sch_S}
  By Theorem \ref{thm_C} \eqref{it_C_A}, we have
   $(x_*, \ep)\in \sC_\sE^\obsm(X)$ for $X\in \Sch_{/S}$ and
   $x\in X(S)$.
}}

\defi{ \label{defi_B0_fai_coh_fun}
Let $X, A\in  \Chp_{/S}$ as before.
\enmt{[\upshape (a)]
\item Let
 $X(A)^\obs\subseteq X(A)$ be a subset.
If for every $x\in X(A)$ such that
 $\ol x_* = (x_*, \ep)\in \sB_0^\obs(X)_0$,
   $x\in X(A)^\obs$,
 then we say that \emph{$X(A)^\obs$ is $\sB_0$-faithful}.

\item
Let $\La\in\Ring_{\Box\tu{-tor}}$.
Recall that we  write $\sD(-)$ for $\sD(-, \La)$.
Let $\Chp_{/S}'\subseteq \Chp_{/S}$ be a subcategory.
A functor  $F: \Chp_{/S}'^\op\ra \Set$ is a
 \emph{cohomological functor} if there exist integers $n_i$ and
 $K_i\in\sD(S)$ and an  inclusions of functors
\eq{ \label{eq_coh_fun}
F \hra \prod_i  H_\et^{n_i}(-, K_i),
}
 where $i$ runs through some (possibly infinite) set.

In particular, we may take $K_i$ to be a commutative $S$-groups $G_i$
 that is $\Box$-torsion.

Note that a cohomological functor is stable.
}}

\lemm{ \label{lemm_B_rec}
Let  $\La\in\Ring_{\Box\tu{-tor}}$.
Let $\obs$ be a cohomological functor
 with the corresponding
 obstruction map $X\mpt X(A)^\obs\subseteq X(A)$
 $($see Section \ref{ptob}$)$.
Then  for every  $X\in \Chp_{/S}^\Box$,
 the subset $X(A)^\obs$ is $\sB_0$-faithful.
}
\pf{
Let $x\in X(A)$ such that
 $\ol x_* = (x_*, \ep)\in \sB_0^\obs(X)_0$.
By definition of  $\sB_0^\obs(X)$
 and Remark \ref{rk_A_struc} \eqref{it_0_1_ess},
 there exists some $x_0\in X(A)^\obs$ and
\eqn{
(x_*, \ep, \mu_*)\ra ({x_0}_*, \ep_0, \ep_0^{-1})
}
 in $\sH_1(X)_1$ with $V((x_*, \ep, \mu_*)) = \ol x_*$,
 where  $\ep_0: q_*\xra{\sim}{x_0}_*p_*$ is
  the natural equivalence.

First we assume
  $\obs$ is of the form
    $\prod_i  H_\et^{n_i}(-, K_i)$.
Since both ${x_0}_*$ and $x_*$ have left adjoints,
 for every $i$,
  the image of $x_0^*: H_\et^{n_i}(X, K_i)\ra H_\et^{n_i}(A, K_i)$
 comes from  $H_\et^{n_i}(S, K_i)$,
 Remark \ref{rk_coh_id} tells us that
 the image of $x^*: H_\et^{n_i}(X, K_i)\ra H_\et^{n_i}(A, K_i)$
 also comes from  $H_\et^{n_i}(S, K_i)$.
For general $\obs$,  by assumption
   it is a subfunctor of  $\prod_i  H_\et^{n_i}(-, K_i)$,
   the previous argument is still correct.
This shows that $x\in X(A)^\obs$.
The proof is complete.
}

Now we apply previous results to local-global obstructions.
Thus let $q: A\ra S$ be $\Spec \bfA_k\ra \Spec k$  induced by  the
 inclusion $k\subset \bfA_k$, where $k$ is a global field.
We have by definition
\ga{ \label{eq_tobs_k}
\XA^{\tobs_\sE} = \Xk\ \bigcup\
 \{x\in \XA \mid \text{ the object
  $\ol x_* = (x_*, \ep)\in \sH_0(X)_0$ } \\ \nonumber
   \text{ lies in the  image of $\sC_\sE^\obsm(X)\hra \sH_0(X)$} \}
}

\thm{ \label{thm_tobs}
Let  $\La\in\Ring_{\Box\tu{-tor}}$.
Let $\Sch_{/k}'\subseteq \Sch_{/k}$ be a full subcategory.
Let $\Esp_{/k}'$ $($resp. $\Chp_{/k}')$ be the corresponding
 full subcategory $($resp. sub $2$-category$)$ $($see
 Remark \ref{rk_Sch_sub}$)$.
Let  $\La\in\Ring_{\Box\tu{-tor}}$.
Let $\sE$ be a  family of maps of $\Chp_{/k}'$  such that
\enmt{[\upshape (a)]
\item $\sE$ contains every degenerate edge,
\item $\sE$  are stable under composition and pullback,
\item every $f\in \sE$ is smooth surjective, and
\item $\sE\subseteq \rep\cap \Box$.
}
Let
\eqn{
\obs: \Sch_{/k}'\ra \Set, \quad X\mpt \XA^\obs
}
 be a map.
To
 every pair $(\obs, \sE)$
 we    associate to
 a subset $X(A)^{\tobs_\sE}$ defined by \eqref{eq_tobs_k},
 produced by $\sC_\sE^\obsm(X)$.
Then $\tobs_\sE$
 defines an obstruction on ${\Chp_{/k}'}^\sE$
 satisfying the following properties.
\enmt{[\upshape (1)]
\item \label{it_tobs_A}
For any $X\in \Sch_{/k}'$,
 we have $\XA^\obs\subseteq \XA^{\tobs_\sE}$.

\item \label{it_tobs_func}
We have  ${\tobs_\sE}$ is functorial on $({\Chp_{/k}'}^\sE)_\rep$
 $($Definition \ref{defi_obs_functorial}$)$.

\item \label{it_tobs_cmp1}
Suppose that we are given two maps
\eqn{
\obs_1\subseteq \obs_2: \Sch_{/k}'\ra \Set
}
Then for every $X\in {\Chp_{/k}'}^\sE$,
 we have $\XA^{\wt{\obs_1}_\sE}\subseteq \XA^{\wt{\obs_2}_\sE}$.

\item \label{it_tobs_B}
If
 for each $f: Y\ra X$ in $\sE\cap (\Sch_{/S})_1$,
  $f$  induces a surjective map
  $Y(A)^\obs\ra X(A)^\obs$,
 then we have the following properties.
\enmt{[\upshape (i)]
\item \label{it_tobs_only}
Let  $X\in {\Chp_{/k}'}^\sE$.
 Suppose that there is a
 a map
 $f: Y\ra  X$
 in $({\Chp_{/k}'}^\sE)_\rep$ with $Y\in \Sch_{/k}'$,
 such that  $\obs$ is the only one for  $Y$
 $($i.e., $\YA^\obs\neq\emptyset$ implies
 $Y(k)^\obs\neq\emptyset)$,
 then  $\tobs_\sE$ is
 the only one for  $X$.

In particular, if $\obs$ is the only one for every $X\in \Sch_{/k}'$,
 then  $\tobs_\sE$ is
 the only one for every $X\in {\Chp_{/k}'}^\sE$.

%
\item \label{it_tobs_Sch}
If moreover,  $\obs$ comes from a cohomological functor
    $\obs: \Chp_{/k}'^\op\ra \Set$,
 then
 for every $X\in {\Chp_{/k}'}^\sE$,
 we have $\XA^{\tobs_\sE}\subseteq \XA^\obs$.

In particular,
 for any $X\in \Sch_{/k}'$,
 we have $\XA^\obs= \XA^{\tobs_\sE}$.
}

\item \label{it_tobs_E1}
Let $\sE_0\subset\sE$ be another family satisfying conditions
 at the beginning of the theorem.
Then for any $X\in {\Chp_{/k}'}^{\sE_0}$,
 we have $\XA^{\tobs_{\sE_0}}\subseteq \XA^{\tobs_\sE}$.

\item \label{it_tobs_Fobs_Chp}
Suppose that
 $\obs$ extends  to a map
\eqn{
  \obs: \Chp_{/k}'^\sE\ra \Set, \quad X\mpt X(A)^\obs\subseteq X(A).
}
Let  $X\in {\Chp_{/k}'}^\sE$ such that there is
 a map
 $f: Y\ra  X$
 in $({\Chp_{/k}'}^\sE)_\rep$ with $Y\in \Sch_{/k}'$,
  and   $f$  induces a surjective map
  $Y(A)^\obs\ra X(A)^\obs$,
Then   we have $\XA^\obs \subseteq\XA^{\tobs_\sE}$.

In particular, if for
  each $f: Y\ra X$ in $({\Chp_{/k}'}^\sE)_\sE$ with $Y \in \Sch_{/k}'$,
  $f$  induces a surjective map
  $Y(A)^\obs\ra X(A)^\obs$,
 then
 for every $X\in {\Chp_{/k}'}^\sE$,
 we have $\XA^\obs \subseteq\XA^{\tobs_\sE}$.
If moreover $\obs$ comes from a cohomological functor
    $\obs: \Chp_{/k}'^\op\ra \Set$,
 then
 for every $X\in {\Chp_{/k}'}^\sE$,
 we have $\XA^\obs= \XA^{\tobs_\sE}$.

\item \label{it_tobs_ttobs}
We have an equivalence  of functors
\eqn{
\sC_\sE^\obsm\xra{\sim} \sC_\sE^{\tobs_\sE,m}:
 \N({\Chp_{/k}'}^\sE)_\rep^\op\ra \Cat_\infty.
}
In particular, for every $X\in {\Chp_{/k}'}^\sE$,
 we have $\XA^{\wt{(\tobs_\sE)}_\sE} = \XA^{\tobs_\sE}$.
}}
\pf{
The proof generally follows from Theorem \ref{thm_C}
 along with Remark \ref{rk_Sch_sub},
 and the definition
 \eqref{eq_tobs_k}.

The property \eqref{it_tobs_A} obviously follows from
 Theorem \ref{thm_C} \eqref{it_C_A}.

For \eqref{it_tobs_func},  by Theorem \ref{thm_C} \eqref{it_C_sub1},
 we have inclusion of functors
\eq{ \label{eq_sub1}
\sC_\sE^\obsm \subseteq \sH_0:
 \N(\Chp_{/k}^\sE)_\rep^\op\ra \Cat_\infty.
}
Since the functor $\sC_\sE^\obsm$ on edges
 agrees with $*$-pullback, with right adjoint agrees with $*$-pushforward,
 the functoriality of $\tobs_\sE$ then follows from definition.

The property \eqref{it_tobs_cmp1} follows directly from
 Theorem \ref{thm_C} \eqref{it_C_cmp1}.

To show \eqref{it_tobs_only},
 let $X\in  {\Chp_{/k}'}^\sE$ such that $\XA^{\tobs_\sE}\neq\emptyset$.
We want to show that  $\Xk\neq\emptyset$.
Suppose that  $\Xk=\emptyset$.
Then by definition \eqref{eq_tobs_k}, $\sC_\sE^\obsm(X)\neq\emptyset$.
Again by Theorem \ref{thm_C} \eqref{it_C_sub1},
 the functor $\sC_\sE^\obsm$ on edges
 agrees with $*$-pullback, with right adjoint agrees with $*$-pushforward.
Thus  we have $\sC_\sE^\obsm(Y)\neq \emptyset$.
By the argument in the proof of
  Theorem \ref{thm_C} \eqref{it_C_B},
    we have a fully faithful  functor
 $\sC_\sE^\obsm(Y)\hra \sB_0^\obs(Y)$,
 it follows that
 $\sB_0^\obs(Y)\neq \emptyset$.
Then by definition of $\sB_0^\obs$ and $\sA_0^\obs$,
 we have
 $\YA^\obs\neq\emptyset$.
Since by assumption $\obs$  is the only one for $Y$.
 it follows that
 $Y(k)\neq\emptyset$.
Note we have the map $Y\xra{f}  X$.
Thus $\Xk\neq\emptyset$, a contradiction.
For the ``In particular'' part,
 just choose an $\sE$-atlas  for $X$,
 that is  a smooth surjective map
 $f: X_0\ra X_{-1} = X$
 in $({\Chp_{/k}'}^\sE)_{\sE}$ with $X_0\in {\Esp_{/k}'}^\sE$,
 and  choose an $\sE$-atlas  for $X_0$,
 that is  a {\etale} surjective map
 $f_0: X_{0,0}\ra X_{0,-1} = X_0$
 in $({\Esp_{/k}'}^\sE)_{\sE}$ with $X_{0,0}\in \Sch_{/k}'$.
Then the argument is similar.
Thus we have shown
 \eqref{it_tobs_only}.


For \eqref{it_tobs_Sch}, suppose that $x\in\XA^{\tobs_\sE}$.
By definition if $x\in \Xk$, trivially we have $x\in \XA^\obs$.
Otherwise, we may assume that
 $\ol x_* = (x_*, \ep)\in \sC_\sE^\obsm(X)_0$.
Since $\obs$ here is functorial
 on  $\Chp_{/S}^\sE$
 (Remark \ref{rk_functorial}),
 by Theorem \ref{thm_C} \eqref{it_C_B},
 for every $X\in {\Chp_{/k}'}^\sE$,
  we have a fully faithful  functor
 $\sC_\sE^\obsm(X)\hra \sB_0^\obs(X)$.
Thus we may further assume that
 $\ol x_*\in \sB_0^\obs(X)_0$.
Since $\obs$ comes from a cohomological functor,
  by Lemma \ref{lemm_B_rec},
  $\XA^\obs$ is $\sB_0$-faithful.
Thus we also have
 $x\in \XA^\obs$.
This shows that
 $\XA^{\tobs_\sE}\subseteq \XA^\obs$.
Then combining \eqref{it_tobs_A},
 for any $X\in \Sch_{/k}'$, that
 $\XA^\obs= \XA^{\tobs_\sE}$ is clear.

The property \eqref{it_tobs_E1} also follows directly from
 Theorem \ref{thm_C} \eqref{it_C_E1}.

For \eqref{it_tobs_Fobs_Chp},
By the assumption,
 we know that $\YA^\obs\ra \XA^\obs$ is surjective.
By  \eqref{it_tobs_A},
 we have $\YA^\obs\subseteq \YA^{\tobs_\sE}$.
Noting that $\tobs_\sE$ is functorial
 on $({\Chp_{/k}'}^\sE)_\rep$
 (by \eqref{it_tobs_func}),
 it follows that
\eqn{
\XA^\obs=f(\YA^\obs)\subseteq  f(\YA^{\tobs_\sE})
 \subseteq \XA^{\tobs_\sE}.
}
For the ``In particular'' part,
 just choose an $\sE$-atlas
 $f: X_0\ra X_{-1} = X$,
 and then choose an $\sE$-atlas
 $f_0: X_{0,0}\ra X_{0,-1} = X_0$
Then the argument is similar.
The ``If moreover" part follows from
  \eqref{it_tobs_Sch}.

For \eqref{it_tobs_ttobs},
 note that for
  every $X\in \Sch_{/k}'$,
 we have  fully faithful  functors
\eqn{
  \sA^{\tobs_\sE}\hra \sC_\sE^\obsm(X)\hra
  \sC_\sE^{\tobs_\sE,m}(X)\hra
   \sH_0(X),
}
 the first of which is due to the definition of $\tobs_\sE$
 and Remark \ref{rk_tobs} \eqref{it_Sch_S},
 the second due to \eqref{it_tobs_A} and Theorem \ref{thm_C} \eqref{it_C_cmp1},
 and the last due to Theorem \ref{thm_C} \eqref{it_C_sub1}.
Since
\eqn{
\sC_\sE^\obsm\in
 \Fun^{\RAd, \sE}(\N({\Sch_{/k}'})_\rep^\op, \Cat_\infty)_{/\sH_0},
}
 by the construction of $\sC_\sE^{\tobs_\sE,m}$ \eqref{eq_C},
 obtain the equivalence $\sC_\sE^\obsm\xra{\sim} \sC_\sE^{\tobs_\sE,m}$.
The proof is complete.
}

\rk{ \label{rk_tobs_cap_coh}
\enmt{[\upshape (1)]
\item \label{it_tobs_cap}
  If $\XA^\obs = \bigcap_i \XA^{\obs_i}$, we can apply Theorem \ref{thm_tobs}
  with different $\Box_i$ and $\La_i$ but the same $\sE$
  for each $\obs_i$ to obtain $\wt{\obs_i}_\sE$.
  One can check that   Theorem \ref{thm_tobs} is still correct for
   the obstruction $X\mpt \bigcap_i \XA^{\wt{\obs_i}_\sE}$, which
   we still write  $\tobs_\sE$.
\item \label{it_tobs_coh}
In the definition of a cohomological functor \ref{defi_B0_fai_coh_fun},
 we may take $\Box_i$ for every $i$ and
 allow $K_i\in \sD(S, \La_i)$  where $\La_i$ is $\Box_i$-torsion.
Since
 \eqn{
   \XA^{\prod_i  H_\et^{n_i}(-, K_i)} = \bigcap_i \XA^{H_\et^{n_i}(-, K_i)},
 }
 by \eqref{it_tobs_cap}, Theorem \ref{thm_tobs}  is still correct.
}}

\section{Concrete examples} \label{conc}

Throughout this section,
 let $q: A\ra S$ be $\Spec \bfA_k\ra \Spec k$  induced by  the
 inclusion $k\subset \bfA_k$, where $k$ is a global field.

\prop{ \label{prop_t_all}
Let $\Sch_{/k}'\subseteq \Sch_{/k}$ be a full subcategory.
Let $\Esp_{/k}'$ $($resp. $\Chp_{/k}')$ be the corresponding
 full subcategory $($resp. sub $2$-category$)$ $($see
 Remark \ref{rk_Sch_sub}$)$.
Let $\sE$ be a  family of maps of $\Chp_{/k}'$  such that
\enmt{[\upshape (a)]
\item $\sE$ contains every degenerate edge,
\item $\sE$  are stable under composition and pullback,
\item every $f\in \sE$ is smooth surjective, and
\item $\sE\subseteq \rep\cap \Box$.
}
Then we have the following corresponding obstructions,
 functorial on $({\Chp_{/k}'}^\sE)_\rep$ and preserving original
 relations,
\gan{
  \XA^{\tBr_\sE} =  \XA^{\tconn_\sE} = \XA^{\tsdesc_\sE} = \XA^{\thoZ_\sE}
   \supseteq \\
 \XA^{\tho_\sE} =  \XA^{\tddesc_\sE} =  \XA^{\tfdesc_\sE} =
   \XA^{\tetBr_\sE} = \XA^{\tdesc_\sE}.
}}
\pf{
Fix   any $\La\in\Ring_{\Box\tu{-tor}}$ and apply
 Theorem \ref{thm_tobs} for $\obs = \Br$, $\conn$, $\sdesc$,
  $h\ZZ$, $h$, $(\ddesc)$,  $(\fdesc)$, $(\etBr)$, $\desc$.
}

Before we give some examples related to $(\etBr)$,
 we need to describe some classes of maps.
\defi{ \label{defi_torsorial}
\enmt{[\upshape (a)]
\item A \emph{torsor tree} is a tree with each vertex being an object $X\in
  \Chp_{/k}$, and for each vertex $X$, all its children are of the form
  $Y^\s\xra{G^\s} X$, $\s\in H^1(k, G)$, where  $G$ is a smooth
  linear $k$-group
  and $Y\xra{G} X$ is  a $G$-torsor.
  Note that we allow different $G$ in one torsor tree.
\item A map $f: Y\ra X$ in  $\Chp_{/k}$ is \emph{torsorial} if  there is a
  torsor tree with $X$ being equivalent to
  its root, and $f: Y\ra X$ being equivalent to the coproduct of
  all its leaves.
\item  A map $f: Y\ra X$ in  $\Chp_{/k}$ is \emph{quasi-torsorial} if $f$ is
  smooth, representable, and there
  is some torsorial map $Z\ra X$ that factorizes through $f$.
}}

We note that $(\etBr)$ is defined on $\Chp_{/k}$ (see Chapter \ref{ptob}).
\lemm{ \label{lemm_quasi_tor}
We have the following statements.
\enmt{[\upshape (1)]
\item \label{it_qt_com_bc}
  Quasi-torsorial maps are stable under composition and pullback.
\item \label{it_qt_etbr_surj}
 Let $f: Y\ra X$ a quasi-torsorial map
  such that $X$ is a smooth $k$-variety.
 Then  $f$  induces a surjective map
  $Y(A)^\etBr\ra X(A)^\etBr$.
\item \label{it_qt_srr}
 Every torsorial map is quasi-torsorial, and
 every quasi-torsorial map is smooth, surjective, and representable.
}}
\pf{
  For \eqref{it_qt_com_bc},
  let $f: Y\ra X$  and $g: Z\ra Y$ be two quasi-torsorial maps.
  By definition there are maps $u: Y'\ra Y$ and $v: Z'\ra Z$ such that
   $f'=fu$ and $g'=gv$ are torsorial maps, as depicted in the following
   $2$-commutative diagram in $\Chp_{/k}$
  \eqn{\xymatrix{
  Z''\ar[r]^-{g''}\ar[d]^-{w} &Y'\ar[r]^-{f'}\ar[d]^-{u} &X \\
  Z'\ar[r]^-{g'}\ar[d]^-v &Y\ar[ur]_-f \\
  Z\ar[ur]_-g
  }}
   with the square being $2$-{\Cart}.
  One  checks that $f'g''$ is torsorial, which
   shows that $fg$ is quasi-torsorial, using the fact that coproducts
   commutes with  pullbacks in $\Chp$.
  The verification about stability of  pullback  is left to the reader.

  For \eqref{it_qt_etbr_surj}, we first show the result for torsorial map.
  Let  $f: Y\ra X$ be  a torsorial map and $X$  a smooth $k$-variety.
  We will see that $Y$ is at least a smooth $k$-scheme.
  Note that in general $Y$ is not  a variety since it not necessary
  quasi-compact.
  Nevertheless,
   by definition, there is a torsor tree with $X$ being equivalent to its root,
  and $f: Y\ra X$ being equivalent to the coproduct of  all its leaves.
  We will show that every vertex of the tree is a smooth $k$-variety and
   that   $f$  induces a surjective map
   $Y(A)^\etBr\ra X(A)^\etBr$.
  By the  fact that
   $(\etBr)$ is functorial, it suffices to show that
   every vertex of the tree is a smooth $k$-variety and
   every $x\in X(A)^\etBr$
   lifts to $Y_0(A)^\etBr$ for some leaf $Y_0$.
  By induction on the  depth of the tree
   and  the problem reduces to
  the case of depth one.
  Thus we may assume that  $f$ is of the form
  \eqn{
    \coprod_{\s\in H^1(k, G)} f^\s: \coprod_{\s\in H^1(k, G)} Y^\s\ra
     X,
  }
   where where  $G$ is a smooth linear $k$-group
  and $Y\xra{G} X$ is  a $G$-torsor.
  It is clear that $Y^\s$ is a smooth $k$-variety for all $\s\in H^1(k, G)$.
  Moreover, by \cite[Thm. 1.1]{cao20sous},
   the induced map $\coprod_\s Y^\s(A)^\etBr\ra X(A)^\etBr$ is surjective.
  This  shows the result for torsorial map.
  For for quasi-torsorial map, the result follows from the torsorial case and
   the functoriality of $(\etBr)$.

  For \eqref{it_qt_srr},  we first  show      that  every
   torsorial map is smooth, surjective and representable.
  Consider the pullback $f$ of a torsorial map   by a scheme.
  As in the argument of \eqref{it_qt_etbr_surj}, we may assume that
   $f$ is of the form
  \eqn{
    \coprod_{\s\in H^1(k, G)} f^\s: \coprod_{\s\in H^1(k, G)} Y^\s\ra
     X,
  }
   where where  $X$ is a scheme, $G$ is a smooth linear $k$-group
   and $Y\xra{G} X$ is  a $G$-torsor.
  Since for every $\s$, $Y^\s\ra X$ is  a $G^\s$-torsor,
   in particular it is smooth, surjective  map between schemes.
  This shows that a torsorial map is surjective and quasi-torsorial.
  Next, by definition a quasi-torsorial map is
   smooth, representable, and   clearly is also surjective since
   a torsorial one is surjective.

  The proof is complete.
}

\prop{ \label{prop_t_etbr}
Let $k$ be a number field.
Let $\Sch_{/k}'\subseteq \Sch_{/k}$ be
 the full subcategory spanned by
 smooth $k$-varieties.
Let $\Esp_{/k}'$ $($resp. $\Chp_{/k}')$ be the corresponding
 full subcategory $($resp. sub $2$-category$)$ $($see
 Remark \ref{rk_Sch_sub}$)$.
Let $\sE$ be the  family of all quasi-torsorial maps $($see Definition
 \ref{defi_torsorial}$)$.
Then $\tetBr_\sE$ is defined on  ${\Chp_{/k}'}^\sE$ and
 functorial on $({\Chp_{/k}'}^\sE)_\rep$, with the following properties.

\enmt{[\upshape (1)]

\item \label{it_t_etbr_only}
Let  $X\in {\Chp_{/k}'}^\sE$.
Suppose that there is
 a map
 $f: Y\ra  X$
 in ${\Chp_{/k}'}^\sE$ with $Y\in \Sch_{/k}'$,
 such that  $(\etBr)$ is the only one for  $Y$,
 then  $\tetBr_\sE$ is
 the only one for  $X$.

\item \label{it_t_etbr_B}
If we shrink $\Sch_{/k}'\subseteq \Sch_{/k}$ to
 the full subcategory spanned by
 torsors over $k$ under linear connected $k$-groups,
 then for every $X\in {\Chp_{/k}'}^\sE$,
 $\tetBr_\sE$ is the only one for  $X$ and
 we have $\XA^{\tetBr_\sE}\subseteq \XA^{\Br_\tor}$,
 where  $\Br_\tor$ is the torsion part of $\Br$
  (which coincides with   $\XA^{\Br}$ if $X\in \Sch_{/k}'$).
In particular,
 for any $X\in \Sch_{/k}'$,
 we have
\gan{
  \XA_{\Br_\tor} = \XA^{\Br} =  \XA^{\Br_1}= \\
  \XA^{\Br_{2/3}} =
  \XA^\desc =  \XA^\etBr= \XA^{\tetBr_\sE},
}
 where $\Br_{2/3}$  is defined by \cite[(2.2)]{cao24sous}.
}}
\pf{
  We need to verify conditions for $\sE$ in Theorem \ref{thm_tobs}.
  By taking $G=1$ in the torsor tree, we see that (a) is clear.
  Condition (b)  follows from  Lemma  \ref{lemm_quasi_tor}
   \eqref{it_qt_com_bc}.
  Conditions (c) and (d) just follows from the definition.
  We set $\Box=\{\ell\}$ where $\ell$ is  any fixed prime
  and take $\La=\ZZ/\ell\ZZ$.
  Since  $\ch k = 0$,
   every object in $\Chp_{/k}$ is $\Box$-coprime,   $\La$ is
   $\Box$-torsion and  $\sE\subseteq     \Box$.
  Thus $\sE$ works and we can  apply Theorem \ref{thm_tobs}  with
   the previous $\sE$, $\Box$, $\La$ and $\obs=(\etBr)$.
  Then $\tetBr_\sE$ is defined on  ${\Chp_{/k}'}^\sE$ and
    the functoriality on $({\Chp_{/k}'}^\sE)_\rep$ follows from
   Theorem \ref{thm_tobs}   \eqref{it_tobs_func}.

  For \eqref{it_t_etbr_only},
   we see that  the  condition of Theorem \ref{thm_tobs}
   \eqref{it_tobs_B} holds by
   Lemma  \ref{lemm_quasi_tor} \eqref{it_qt_etbr_surj}.
  Then the results follows from
   Theorem \ref{thm_tobs}   \eqref{it_tobs_only}.

  For \eqref{it_t_etbr_B},  first note that for any $X\in \Sch_{/k}'$,
   $\Br$ is the  only one for $X$ \cite[Thm. 2.2]{borovoi96bm}.
   Thus $(\etBr)$ is also the only one for $X$ since
    $\XA^\etBr\subseteq\XA^{\Br}$.
  We also have
   \cite[Cor. 1.4]{cao24sous}
  \gan{
    \XA^{H_\et^2(-, \mu_\infty)} =
     \XA^{\Br_\tor} =\XA^{\Br} =  \XA^{\Br_1} = \\
     \XA^{\Br_{2/3}} =
    \XA^\desc =  \XA^\etBr,
  }
   where
  \eqn{
    H_\et^2(-, \mu_\infty) = \prod_n H_\et^2(-, \mu_n).
  }
  It follows that
  \eqn{
    \tetBr_\sE = \wt{\Br_\tor}_\sE = \wt{  H_\et^2(-, \mu_\infty) }_\sE
  }
    on  ${\Chp_{/k}'}^\sE$ and is the only one by
    ``In particular" part of Theorem \ref{thm_tobs}   \eqref{it_tobs_only}.
  Since   $H_\et^2(-, \mu_\infty)$ is cohomological on $\Chp_{/k}$ if we
   take $\La_n = \ZZ/n\ZZ$   is cohomological on $\Chp_{/k}$
    (see Remark \ref{rk_tobs_cap_coh}
       \eqref{it_tobs_coh}),
   we use  Theorem \ref{thm_tobs}
   \eqref{it_tobs_Sch} to deduce
  \eqn{
    \XA^{\tetBr_\sE} =  \wt{  H_\et^2(-, \mu_\infty) }_\sE
    \subseteq
    \XA^{H_\et^2(-, \mu_\infty)} = \XA^{\Br_\tor}.
  }
  The proof is complete.
}

\lemm{ \label{lemm_contracted_rep}
Let $S$ be a scheme, $Y$  a quasi-affine  $S$-schemes,
 $G$ a smooth affine group  scheme over $S$ that acts on $Y$,
 $Z\xra{G} B$ a $G$-torsor in $S$-schemes.
Then the contracted product $Z\tm_S^G Y$  is
 represented by a quasi-affine $B$-scheme.

Moreover, if $Y\ra S$ is
 separated $($resp. smooth, resp. of finite type, resp.
 geometrically integral$)$,
 then $Z\tm_S^G B\ra Y$  is also
 separated $($resp. smooth, resp. of finite type, reps.
 geometrically integral$)$.
}
\pf{
  The action $Z\tm_S G\ra Z$ sending $(g,z)$ to $zg^{-1}$ makes
   $Z\xra{G} B$ a right $G$-torsor.
  Thus there exist an fppf map of schemes $B'\ra B$ trivializing $Z$,
   that is, $Z\tm_B B' \xra{\sim} G\tm_B B'$
   with right $G\tm_B B'$-action.
  By definition, $Z\tm_S^G Y$ is the quotient stack
   $[(Z\tm_S Y)/G]$ where $G$ acts on $Z\tm_S Y$ via the
   diagonal $g(z,y) = (zg^{-1}, gy)$.
  This action is free since the $G$-action on $Z$ is.
  It follows from \cite[Cor. 4.6.8 (2)]{moduli} that
   $[(Z\tm_S Y)/G]$ is represented
   by the quotient space $X=(Z\tm_S Y)/G$.

  Consider the projection $X\ra Z/G\cong B$.
  As fppf-sheaves, we shall show that it is a fiber bundle with fiber $Y$.
  Indeed,
  \gan{
  X \tm_B B'\xra{\sim} ((Z\tm_S Y)\tm_B B' )/G\xra{\sim}
   ((Z\tm_B B')\tm_S Y)/G \xra{\sim} \\
   (G\tm_S B'\tm_S Y)/G \xra{\sim} (G\tm_S ( Y\tm_S B'))/G
   \xra{\sim} Y\tm_S B'
  }
   compatible with the projection to $B'$, making $X\ra B$ a fiber
   bundle.
  Note that the right hand side is a quasi-affine $B'$-scheme.
  One can verify that this also  equip the $B'$-scheme $Y\tm_S B'$
   with a descent datum under $B'\ra B$.
  It follows by {\Grot}'s fppf descent
   (see, e.g., \cite[Thm. 4.3.5 (2)]{poonen17rational})
   that there is a quasi-affine $B$ scheme $X_0$,
   such that  $X_0\tm_B B'\xra{\sim} Y\tm_S B'$.

  The proof of the first statement
   is complete if we show that $X\cong X_0$
   as fppf $B$-sheaves, which follows from the
   fact that the functor $\Sh: (\tu{Schemes}_{/S})_\fppf
    \ra (\tu{Groupoids})$ is a stack
    \cite[Thm. 4.2.12]{olsson16stack}.
  The  second statement is correct by
   descent properties of fppf maps in schemes,
   see, for example, \cite[Thm. 4.3.7]{poonen17rational}.
}
\rk{
 Lemma \ref{lemm_contracted_rep}  is a generalization of
  \cite[Lem. 2.2.3]{torsor} (or \cite[6.5.6.3]{poonen17rational})
  which it is the special case where $Y$ is an affine  $S$-scheme
  and $B=S$.
}

\cor{ \label{cor_t_etbr_only}
Let $Y$ be a smooth quasi-affine variety over a number field $k$,
 $G$ be a linear $k$-group with an action on $Y$,
 and $X = [Y/G]$ the quotient stack.
Since $G$ is linear, we may embed $G$ into  $\SL_n$  for some $n$.
Then the contracted product $Y' = \SL_n\tm_k^G Y$  is a smooth variety.

Let $\sE$ be the  family of all quasi-torsorial maps.
Then $\tetBr_\sE$ is defined on $X$.
If moreover, $(\etBr)$ is the only obstruction for $Y'$,
 then $\tetBr_\sE$ is the only one for $X$.
}
\pf{
Applying Lemma \ref{lemm_contracted_rep} with $S = \Spec k$,
 $Z = \SL_n$ and $B= \SL_n/G$,
 we obtain that  $Y' = \SL_n\tm_k^G Y$  is a smooth variety.

But we have an equivalence $X = [Y/G]\xra{\sim} [Y'/\SL_n]$ since
 $Y'\xra{\SL_n} X$ is an $\SL_n$-torsor
 (or, c.f. \cite[Ex. 3.4.19]{moduli}).
Under notation of Proposition \ref{prop_t_etbr},
 since $H^1(k, \SL_n)=\{*\}$, it follows that  $Y'\ra X$ is
 an $\sE$-atlas for $X$ with $Y'\in \Sch'_{/k}$.
It follows that $X\in {\Chp_{/k}'}^\sE$.
The remaining assertions  follows from
 Proposition \ref{prop_t_etbr},
  \eqref{it_t_etbr_only}.
}


\section*{Acknowledgment}
The author would like to thank
 Weizhe Zheng,
 Yang Cao,
 Junchao Shentu and
 Han Wu
 for helpful discussions.
The work was partially done during the author's visit to
 the Morningside Center of Mathematics, Chinese
 Academy of Sciences. He thank the Center for its hospitality.

\bibliography{unibib}
\bibliographystyle{amsalpha}
\end{document}